\documentstyle[11pt]{article}
\def\Bbb R{{\rm \bf R}}
\def\proclaim#1{\vskip2mm{\bf #1}\em}
\def\endproclaim{\em \vskip2mm}
\def\tag#1{\eqno(#1)}
\def\gathered{\begin{array}{c}}
\def\endgathered{\end{array}}
\def\text{\mbox}

\begin{document}

\title {Reconstruction of a source domain from the Cauchy data: II. Three dimensional case}
\author{Masaru IKEHATA\footnote{
Laboratory of Mathematics,
Graduate School of Advanced Science and Engineering,
Hiroshima University,
Higashihiroshima 739-8527, JAPAN}
\footnote{Emeritus Professor at Gunma University}
}
\maketitle

\begin{abstract}
This paper is concerned with reconstruction issue of some typical inverse problems and consists of three parts.  First a framework of the enclosure method
for an inverse source problem governed by the Helmholtz equation
at a fixed wave number in three dimensions is introduced.  It is based on the nonvanishing of the coefficient of the leading profile
of an oscillatory integral over a domain having a conical singularity. 
Second an explicit formula of the coefficient for a domain having a circular cone singularity 
and its implication under the framework are given. 
Third, an application under the framework to an inverse obstacle problem 
governed by an inhomogeneous Helmholtz equation at a fixed wave number in three dimensions is given.

\noindent
AMS: 35R30

$\quad$



\end{abstract}


\section{Introduction}

More than twenty years ago, in \cite{Ik} the author obtained the extraction formula of the support function 
of an unknown polygonal source domain in an inverse source problem governed by the Helmholtz equation
and polygonal penetrable obstacle in an inverse obstacle problem governed by an inhomogeneous Helmholtz eqution.
All the problems considered therein are in two dimensions and employ only a single set of Cauchy data
of a solution of the governing equation at a fixed wave number in a bounded domain.
Those results can be considered as the first application of a single measurment version of 
the {\it enclosure method} introduced in \cite{Ik0}.

Succeding to \cite{Ik}, in \cite{IkC} the author found another unexpected application
of the enclosure method out to the Cauchy problem for the stationary Schr\"odinger equation
$$\displaystyle
-\Delta u+V(x)u=0
\tag {1.1}
$$
in a bounded domain $\Omega$ of $\Bbb R^n$, $n=2,3$.  Here $V\in L^{\infty}(\Omega)$ and both $u$ and $V$ can be complex valued functions.
We established an explicit representation or computation formula for an arbitrary solution $u\in H^2(\Omega)$ to the equation (1.1) in $\Omega$
in terms of its Cauchy data on a part of $\partial\Omega$.  See also \cite{IS} for its numerical implementation.
Note also that the idea in \cite{IkC} has been applied
to an inverse source problem governed by the heat equation together with an inverse heat conduction problem
in \cite{IkH}, \cite{IkHC}, respectively.

The idea introduced therein is to make use of the complex geometrical optics solutions (CGO) with a large parameter $\tau$ for the modified equation
instead of (1.1):
$$\begin{array}{ll}
\displaystyle
-\Delta v+V(x)v=\chi_{D_y}(x)v, & x\in\Omega,
\end{array}
$$
where $y$ is a given point in $\Omega$,  $D_y\subset\subset\Omega$ is the inside of a triangle, tetrahedron for $n=2,3$, respectively 
with a vertex at $y$ and $\chi_{D_y}(x)$ is the characteristic function of $D_y$.
The solution is the same type as constructed one in \cite{SU1} for $n=2$, \cite{SU2} for $n=3$ and has the following form as $\tau\rightarrow\infty$
$$\displaystyle
v\sim e^{x\cdot z},
$$
where $z=\tau(\omega+i\vartheta)$ and both $\omega$ and $\vartheta$ are unit vectors perpendicular to each other.
This right-hand side is just the complex plane wave used in the Calder\'on method \cite{C}.

Note that, in \cite{IkS} another simpler idea to make use of the CGO solutions of another modified equation described below is
presented:
$$\begin{array}{ll}
\displaystyle
-\Delta v+V(x)v=\chi_D(x)e^{x\cdot z}, & x\in\Omega.
\end{array}
$$
Using integration by parts we reduced the problem of computing the value of $u$ at given point $y$, essentially, to clarifying the leading profile of the following oscillatory integral
as $\tau\rightarrow\infty$:
$$\displaystyle
\int_{D_y} e^{x\cdot z}\,\rho(x)dx,
$$
where $\rho(x)$ is uniformly H\"older continuous on $\overline{D_y}$\footnote{In this case $\rho=u$.}.
Note that the asymptotic behaviour of this type of oscillatory integral in {\it two dimensions} is the key point of the enclosure method developed in \cite{Ik}.

In \cite{IkC} we clarified the leading profile in more general setting as follows.
Given a pair $(p,\omega)\in\Bbb R^n\times S^{n-1}$ and $\delta>0$ let $Q$ be an arbitrary non empty bounded open subset
of the plane $x\cdot\omega=p\cdot\omega-\delta$ with respect to the relative topology from $\Bbb R^n$.
Define the bounded open subset of $\Bbb R^n$ by the formula
$$\displaystyle
D_{(p,\omega)}(\delta,Q)
=\cup_{0<s<\delta}\,
\left\{p+\frac{s}{\delta}(z-p)\,\left\vert\right.\,z\in Q\,\right\}.
\tag {1.2}
$$
This is a cone with the base $Q$ and apex $p$, and lying in the slab $\{x\in\Bbb R^n\,\vert\,p\cdot\omega-\delta<x\cdot\omega<p\cdot\omega\,\}$.
Note that $\delta=\text{dist}\,(\{p\},Q)$ is called the height.  If $Q$ is given by the inside of a polygon, the cone (1.2) is called a {\it solid pyramid}.
In particular, if $Q$ is given by the inside of a triangle, cone (1.2) becomes a tetrahedron.

On (2.2) in \cite{IkC} we introduced a special complex constant associated with the domain (1.2) which is given by
$$\displaystyle
C_{(p,\omega)}(\delta, Q,\vartheta)=2s\int_{Q_s}\frac{dS_z}{\{s-i(z-p)\cdot\vartheta\}^n},
\tag {1.3}
$$
where $i=\sqrt{-1}$, $0<s<\delta$ and $Q_s=D_{(p,\omega)}(\delta,Q)\cap\{x\in\Bbb R^n\,\vert x\cdot\omega=p\cdot\omega-s\,\}$
and the direction $\vartheta\in S^{n-1}$ is perpendicular to $\omega$.
Note that in \cite{IkC} complex constant $C_{(p,\omega)}(\delta,Q,\vartheta)$ is simply written as $C_D(\omega,\omega^{\perp})$ with $\omega^{\perp}=\vartheta$.
As pointed therein out this quantity is independent of the choice $s\in\,]0,\,\delta[$ because of
the one-to-one correspondence between $z\in Q_s$ and $z'\in Q_{s'}$ by the formula
$$
\left\{
\begin{array}{l}
\displaystyle
z'=p+\frac{s'}{s}\,(z-p),
\\
\\
\displaystyle
\displaystyle
dS_{z'}=(\frac{s'}{s})^{n-1}\,dS_z.
\end{array}
\right.
$$

The following lemma describes the relationship between complex constant $C_{(p,\omega)}(\delta, Q,\vartheta)$ and an integral over (1.2).

\proclaim{\noindent Proposition 1.1 (Lemma 2 in \cite{IkC}).}  Let $n=2, 3$.
Let $D=D_{(p,\omega)}(\delta,Q)$ and $\rho\in C^{0,\alpha}(\overline D)$ with $0<\alpha\le 1$.
It holds that, for all $\tau>0$
$$\begin{array}{l}
\displaystyle
\,\,\,\,\,\,
\left\vert
e^{-\tau p\cdot(\omega+i\vartheta)}
\int_D\rho(x)e^{\tau x\cdot(\omega+i\vartheta)}\,dx
-\frac{n-1}{2\tau^n}
\rho(p)\,C_{(p,\omega)}(\delta,Q,\vartheta)\right\vert
\\
\\
\displaystyle
\le\vert\rho(p)\vert\frac{\vert Q\vert}{\delta^{n-1}}
\{(\tau\delta+1)^{n-1}+n-2\}
\frac{e^{-\tau\delta}}{\tau^n}
+\Vert\rho\Vert_{C^{0,\alpha}(\overline D)}
\frac{\vert Q\vert}{\delta^{n-1}}
(\frac{\text{diam}\,D}{\delta})^{\alpha}\frac{C_{n,\alpha}}{\tau^{n+\alpha}},
\end{array}
$$
where $\Vert\rho\Vert_{C^{0,\alpha}(\overline D)}=
\sup_{x,y\in\overline D, x\not=y}\frac{\vert\rho(x)-\rho(y)\vert}{\vert x-y\vert^{\alpha}}$
and
$$\displaystyle
C_{n,\alpha}
=\int_0^{\infty}s^{n-1+\alpha} e^{-s}ds.
$$
\endproclaim

Thus we have, as $\tau\rightarrow\infty$
$$\displaystyle
e^{-\tau p\cdot(\omega+i\vartheta)}
\int_{D_{(p,\omega)}(\delta,Q)}\rho(x)e^{\tau x\cdot(\omega+i\vartheta)}\,dx
=\frac{n-1}{2\tau^n}
\rho(p)\,C_{(p,\omega)}(\delta,Q,\vartheta)+O(\tau^{-(n+\alpha)}).
$$
This is the meaning of complex constant $C_{(p,\omega)}(\delta,Q,\vartheta)$.
Note that the remainder estimate $O(\tau^{-(n+\alpha)})$ is uniform with respect to $\vartheta$.
And also as a direct corollary, instead of (1.3) we have another representation of $C_{(p,\omega)}(\delta,Q,\vartheta)$:
$$\displaystyle
C_{(p,\omega)}(\delta,Q,\vartheta)
=\frac{2}{n-1}
\lim_{\tau\longrightarrow\infty}\tau^ne^{-\tau p\cdot(\omega+i\vartheta)}
\int_{D_{(p,\omega)}(\delta,Q)} e^{\tau x\cdot(\omega+i\vartheta)}dx.
\tag {1.4}
$$
The convergence is uniform with respect to $\vartheta$.

Proposition 1.1 is the one of two key points in \cite{IkC} and gives the role of the H\"older continuity of $\rho$.
Another one is the {\it non-vanishing} of $C_{(p,\omega)}(\delta,Q,\vartheta)$ as a part of the leading coefficient 
of the integral in Proposition 1.1 as $\tau\rightarrow\infty$.  
This is not trivial, in particular, in three dimensional case.  For this we have shown therein the following fact.

\proclaim{\noindent Proposition 1.2(Theorem 2 in \cite{IkC}).}

$\bullet$  If $n=2$ and $Q$ is given by the inside of an {\it arbitrary line segment}, then for all $\vartheta$ perpendicular
to $\omega$ we have  $C_{(p,\omega)}(\delta,Q,\vartheta)\not=0$.

$\bullet$  If $n=3$ and $Q$ is given by the inside of an {\it arbitrary triangle}, then for all $\vartheta$ perpendicular
to $\omega$ we have $C_{(p,\omega)}(\delta,Q,\vartheta)\not=0$.

\endproclaim

The nonvanishing of complex constant $C_{(p,\omega)}(\delta,Q,\vartheta)$
in case $n=2$ has been shown in the proof of Lemma 2.1 in \cite{Ik}.  The proof therein employs
a local expression of the corner around apex as a graph of a function on the line $x\cdot\omega=x\cdot p$ and so the proof
by viewing $D_{(p,\omega)}(\delta,Q)$ as a cone in \cite{IkC} is not developed.

Note that, in the survey paper \cite{IkS} on the enclosure method it is pointed out that ``the Helmholtz version''
of Proposition 1.1 is also valid.
That is, roughly speaking, we have
$$\displaystyle
e^{-p\cdot(\tau\omega+i\sqrt{\tau^2+k^2}\,\vartheta)}
\int_{D_{(p,\omega)}(\delta,Q)}\,\rho(x)e^{x\cdot(\tau \omega+i\sqrt{\tau^2+k^2}\,\vartheta)}\,dx
=\frac{n-1}{2\tau^n}
\rho(p)\,C_{(p,\omega)}(\delta,Q,\vartheta)+O(\tau^{-(n+\alpha)})
\tag {1.5}
$$
with the {\it same constant} $C_{(p,\omega)}(\delta,Q,\vartheta)$, where $k\ge 0$.  
See Lemma 3.2 therein.
The proof can be done by using the same argument as that of Proposition 1.1.  
Note that the function $v=e^{x\cdot(\tau \omega+i\sqrt{\tau^2+k^2}\,\vartheta)}$ satisfies
the Helmholtz equation $\Delta v+k^2 v=0$ in $\Bbb R^n$.

\subsection{Role of nonvanishing in an inverse source problem}

As an application of the nonvanishing of the complex constant $C_{(p,\omega)}(\delta, Q,\vartheta)$, we present here its direct application 
to the inverse source problem considered in \cite{Ik}, however, in {\it three dimensions}.

Let $\Omega$ be a bounded domain of $\Bbb R^3$ with $\partial\Omega\in C^2$.  
We denote by $\nu$ the normal unit outward vector field on $\partial\Omega$.
Let $k\ge 0$. 
Let $u\in H^1(\Omega)$ be an arbitrary weak solution of the Helmholtz equation in $\Omega$ at the wave number $k$:
$$\begin{array}{ll}
\displaystyle
\Delta u+k^2 u=F(x), & x\in\Omega,
\end{array}
\tag {1.6}
$$
where $F(x)$ is an unknown source term such that $\text{supp}\,F\subset\Omega$.
Both $u$ and $F$ can be complex-valued functions.  See \cite{Ik} for the meaning of the solution
and the formulation of the Cauchy data on $\partial\Omega$ in the weak sense.

It is well known that, in general, one can not obtain the uniqueness of the source term $F$ itself from the Cauchy data
of $u$ on $\partial\Omega$. 
In fact, given $\varphi\in C^{\infty}_0(\Omega)$ let $G=F+\Delta\varphi+k^2\varphi$.
We have $\text{supp}\,G\subset\Omega$ and the function $\tilde{u}=u+\varphi$ satisfies
$$\begin{array}{ll}
\displaystyle
\Delta\tilde{u}+k^2\tilde{u}=G(x), & x\in\Omega.
\end{array}
$$
Both $u$ and $\tilde{u}$ have the same Cauchy data on $\partial\Omega$.
It should be pointed out that, however, $F$ and $G$ coincides each other modulo $C^{\infty}$.
This means that the singularity of $F$ and $G$ coincides each other.  This suggests a possibility
of extracting some information about 
a singularity of $F$ or its support from the Cauchy data of $u$ on $\partial\Omega$.

As done in \cite{Ik} in two dimensions, we introduce the special form of the unknown source $F$:
$$F(x)=F_{\rho,D}(x)=
\left\{\begin{array}{lr}
\displaystyle
0, & \quad\text{if $x\in\Omega\setminus D$,}\\
\\
\displaystyle
\rho(x), & \quad\text{if $x\in\,D$.}
\end{array}
\right.
\tag {1.7}
$$
Here $D$ is an unknown non empty open subset of $\Omega$ satisfying $\overline D\subset\Omega$
and $\rho\in L^{2}(D)$ also unknown.  We call $D$ the {\it source domain}, however,
we do not assume the connectedness of not only $D$ but also $\Omega\setminus\overline D$.
The $\rho$ is called the strength of the source.

We are interested in the following problem.

$\quad$

{\bf\noindent Problem 1.}
Extract information about a singularity of the source domain $D$ of $F$ having form (1.7) from the Cauchy data
$(u(x), \frac{\partial u}{\partial\nu}(x))$ for all $x\in\partial\Omega$.

$\quad$

\noindent
Note that we are seeking a {\it concrete procedure} of the extraction.

Here we recall the notion of the regularity of a direction introduced in the enclosure method \cite{Ik}.
The function $h_D(\omega)=\sup_{x\in D}\,x\cdot\omega$, $\omega\in S^{2}$ is called the {\it support function} of $D$.
It belongs to $C(S^2,\Bbb R)$ because of the trivial estimae $\vert h_D(\omega_1)-h_D(\omega_2)\vert\le
\sup_{x\in D}\,\vert x\vert\cdot\vert\omega_1-\omega_2\vert$ for all $\omega_1,\omega_2\in S^2$.
Given $\omega\in S^{2}$,
it is easy to see that the set
$$\displaystyle
H_{\omega}(D)\equiv\left\{x\in \overline D\,\left\vert\right. x\cdot\omega=h_D(\omega)\,\right\}
$$
is non empty and contained in $\partial D$.
We say that $\omega$ is {\it regular} with respect to $D$ if the set $H_{\omega}(D)$ consists of only a single point.
We denote the point by $p(\omega)$.

We introduce a concept of a singularity of $D$ in (1.7).

{\bf\noindent Definition 1.1.}  Let $\omega\in S^{2}$ be regular with respect to $D$.
We say that $D$ has a {\it conical singularity} from direction $\omega$ if
there exists a positive number $\delta$, an open set $Q$ of the plane $x\cdot\omega=h_D(\omega)-\delta$ with respect to the relative topology
from $\Bbb R^3$ such that 
$$\displaystyle
D\cap\left\{x\in\Bbb R^3\,\vert\,h_D(\omega)-\delta<x\cdot\omega<h_D(\omega)\,\right\}=D_{(p(\omega),\omega)}(\delta,Q).
$$

Second we introduce a concept of an {\it activity} of the source term.

{\bf\noindent Definition 1.2.}  Given a point $p\in\partial D$
we say that the source $F=F_{\rho,D}$ given by (1.7) is {\it active} at $p$
if there exist an open ball $B_{\eta}(p)$ centered at $p$ with radius $\eta$,
$0<\alpha\le 1$ and a function $\tilde{\rho}\in C^{0,\alpha}(\overline{B_{\eta}(p)})$ 
such that $\rho(x)=\tilde{\rho}(x)$ for almost all $x\in B_{\eta}(p)\cap D$
and $\tilde{\rho}(p)\not=0$.  Note that $\rho$ together with $\tilde{\rho}$ can be a complex-valued function.

Now let $u\in H^1(\Omega)$ satisfies the equation (1.6) in the weak sense with $F=F_{\rho, D}$ given by (1.7).
Given a unit vector $\omega\in S^2$ define $S(\omega)=\{\vartheta\in S^2\,\vert \omega\cdot\vartheta=0\}$.

Using the Cauchy data of $u$ on $\partial\Omega$, we define the indicator function as \cite{Ik}
$$\displaystyle
I_{\omega,\vartheta}(\tau)=\int_{\partial\Omega}
\left(\frac{\partial u}{\partial\nu}v-\frac{\partial v}{\partial\nu} u\right)\,dS,
$$
where $\vartheta\in S(\omega)$ and
$$\displaystyle
v=e^{x\cdot(\tau\omega+i\sqrt{\tau^2+k^2}\vartheta)},\,\,\tau>0.
$$
And also its derivative with respect to $\tau$
$$\displaystyle
I_{\omega,\vartheta}'(\tau)
=\int_{\partial\Omega}\left(\frac{\partial u}{\partial\nu}\,v_{\tau}-\frac{\partial\,v_{\tau}}{\partial\nu} u\right)\,dS,
$$
where
$$\displaystyle
v_{\tau}=\partial_{\tau}v=\left\{x\cdot\left(\omega+i\frac{\tau}{\sqrt{\tau^2+k^2}}\,\vartheta\,\right)\,\right\}\,v.
$$

The following theorem clarifies the role of the complex constant $C_{(p,\omega)}(\delta, Q,\vartheta)$ in the asymptotic behaviour of the indicator function
together with  its derivative as $\tau\rightarrow\infty$.

\proclaim{\noindent Theorem 1.1.}
Let $\omega$ be regular with respect to $D$ and assume that
$D$ has a conical singularity from direction $\omega$.
Then, we have
$$\displaystyle
\tau^3e^{-\tau h_D(\omega)}e^{-i\sqrt{\tau^2+k^2}p(\omega)\cdot\vartheta}I_{\omega,\vartheta}(\tau)=
\tilde{\rho}(p(\omega))\,C_{(p(\omega),\omega)}(\delta,Q,\vartheta)
+O(\tau^{-\alpha})
\tag {1.8}
$$
and
$$\displaystyle
\tau^3e^{-\tau h_D(\omega)}e^{-i\sqrt{\tau^2+k^2}p(\omega)\cdot\vartheta}I_{\omega,\vartheta}'(\tau)=
\tilde{\rho}(p(\omega))(h_D(\omega)+ip(\omega)\cdot\vartheta)\,C_{(p(\omega),\omega)}(\delta,Q,\vartheta)
+O(\tau^{-\alpha}).
\tag {1.9}
$$
The remainder $O(\tau^{-\alpha})$ is uniform with respect to $\vartheta\in S(\omega)$.
\endproclaim
{\it\noindent Proof.}
Integration by parts yields
$$\displaystyle
I_{\omega,\vartheta}(\tau)=\int_D\rho(x)\,v\,dx
$$
and thus
$$\displaystyle
I_{\omega,\vartheta}'(\tau)=\int_D\rho(x)\,v_{\tau}\,dx.
$$
Recalling Definition 1.1, one has the decomposition
$$\displaystyle
D=D_{(p(\omega),\omega)}(\delta,Q)\cup D',
\tag {1.10}
$$
where 
$$
D'=D\setminus D_{(p(\omega),\omega)}(\delta,Q)\subset\left\{x\in\Bbb R^3\,\vert\,x\cdot\omega\le h_D(\omega)-\delta\,\right\}.
\tag {1.11}
$$
Besides, choosing $\delta$ smaller if necessary, one may assume that $D_{(p(\omega),\omega)}(\delta, Q)\subset B_{\eta}(p(\omega))$, where $\eta$
and $B_{\eta}(p(\omega))$ are same as those of Definition 1.2.

Hereafter we set $p=p(\omega)$ for simplicity of description.
According to the decomposition (1.10), we have
the decomposition of both  $I_{\omega,\vartheta}(\tau)$ and $I_{\omega,\vartheta}'(\tau)$ as follows:
$$\begin{array}{l}
\displaystyle
\,\,\,\,\,\,
e^{-\tau h_D(\omega)}
e^{-i\sqrt{\tau^2+k^2}p\cdot\vartheta}I_{\omega,\vartheta}(\tau)
\\
\\
\displaystyle
=e^{-\tau h_D(\omega)}e^{-i\sqrt{\tau^2+k^2}p\cdot\vartheta}
\int_{D_{(p,\omega)}(\delta, Q)
}\tilde{\rho}(x)\,vdx
+e^{-\tau h_D(\omega)}e^{-i\sqrt{\tau^2+k^2}p\cdot\vartheta}\int_{D'}\rho v dx
\end{array}
\tag {1.12}
$$
and
$$\begin{array}{l}
\displaystyle
\,\,\,\,\,\,
e^{-\tau h_D(\omega)}e^{-i\sqrt{\tau^2+k^2}p\cdot\vartheta}I_{\omega,\vartheta}'(\tau)
\\
\\
\displaystyle
=e^{-\tau h_D(\omega)}e^{-i\sqrt{\tau^2+k^2}p\cdot\vartheta}
\int_{D_{(p,\omega)}(\delta, Q)}\tilde{\rho}(x)\, v_{\tau}dx+
e^{-\tau h_D(\omega)}e^{-i\sqrt{\tau^2+k^2}p\cdot\vartheta}\int_{D'}\rho v_{\tau} dx,
\end{array}
\tag {1.13}
$$
where $p=p(\omega)$.

By (1.11), we see that the second terms on the right-hand sides of (1.12) and (1.13) have the common bound
$O(e^{-\tau\delta}\Vert\rho\Vert_{L^{2}(D)})$.
Thus from (1.5) and (1.12) we obtain (1.8) with
the remainder $O(\tau^{-\alpha})$ which is uniform with respect to $\vartheta\in S(\omega)$.

For (1.13) we write
$$\begin{array}{l}
\displaystyle
\,\,\,\,\,\,
\int_{D_{(p,\omega)}(\delta, Q)}\tilde{\rho}(x)\, v_{\tau}dx
\\
\\
\displaystyle
=\int_{D_{(p,\omega)}(\delta, Q)}\tilde{\rho}(x)\,
\left\{x\cdot\left(\omega+i\frac{\tau}{\sqrt{\tau^2+k^2}}\,\vartheta\,\right)\,\right\}\,v\,dx
\\
\\
\displaystyle
=\int_{D_{(p,\omega)}(\delta, Q)}\tilde{\rho}(x)\, x\cdot\omega\,v\,dx
+i\frac{\tau}{\sqrt{\tau^2+k^2}}\int_{C_{(p,\omega)}(\delta, Q)}\tilde{\rho}(x)\, x\cdot\vartheta\,v\,dx.
\end{array}
$$
Thus applying (1.5) to each of the last terms and using (1.13), we obtain (1.9)
with the remainder $O(\tau^{-\alpha})$ which is uniform with respect to $\vartheta\in S(\omega)$.

\noindent
$\Box$

Thus under the same assumptions as Theorem 1.1, for each $\vartheta\in S(\omega)$ one can calculate 
$$\displaystyle
I(\omega,\vartheta)\equiv \tilde{\rho}(p(\omega))\,\,C_{(p(\omega),\omega)}(\delta,Q,\vartheta)
$$
via the formula
$$\displaystyle
I(\omega,\vartheta)
=\lim_{\tau\rightarrow\infty}\tau^3e^{-\tau h_D(\omega)}e^{-i\sqrt{\tau^2+k^2}p(\omega)\cdot\vartheta}
I_{\omega,\vartheta}(\tau)
\tag {1.14}
$$
by using the Cauchy data of $u$ on $\partial\Omega$ if $p(\omega)$ is known.

As a direct corollary of formulae (1.8) and (1.9), we obtain a partial answer to Problem 1 and the starting point of the main purpose in this paper.

\proclaim{\noindent Theorem 1.2.}
Let $\omega$ be regular with respect to $D$.
Assume that $D$ has a conical singularity from direction $\omega$,
$F_{\rho,D}$ is active at $p=p(\omega)$ and that
direction $\vartheta\in S(\omega)$ satisfies the condition
$$\displaystyle
C_{(p(\omega),\,\omega)}(\delta,Q,\vartheta)\not=0.
\tag {1.15}
$$
Then, there exists a positive number $\tau_0$ such that, for all $\tau\ge\tau_0$
$\vert I_{\omega,\vartheta}(\tau)\vert>0$ and we have the following three asymptotic formulae.
The first formula is
$$\displaystyle
\lim_{\tau\longrightarrow\infty}\frac{\log\vert I_{\omega, \vartheta}(\tau)\vert}{\tau}=h_D(\omega)
\tag {1.16}
$$
and second one
$$\displaystyle
\lim_{\tau\rightarrow\infty}
\frac{I_{\omega,\vartheta}'(\tau)}{I_{\omega,\vartheta}(\tau)}
=h_D(\omega)+i\,p(\omega)\cdot\vartheta.
\tag {1.17}
$$
The third one is the so-called $0$-$\infty$ criterion:
$$\displaystyle
\lim_{\tau\longrightarrow\infty}e^{-\tau t}\vert I_{\omega, \vartheta}(\tau)
\vert
=
\left\{
\begin{array}{ll}
0, & \text{if $t\ge h_D(\omega)$,}\\
\\
\displaystyle
\infty, & \text{if $t<h_D(\omega)$.}
\end{array}
\right.
\tag {1.18}
$$

\endproclaim

This provides us the framework of the approach using the enclosure method for the source domain with a conical singularity
from a direction.

Some remarks are in order.

$\bullet$  In two dimensions, by Proposition 1.2 the condition (1.15) is redundant
and we have the same conclusion as Theorem 1.2.

$\bullet$  The formula (1.17) is an application of the idea 
``taking the logarithmic derivative of the indicator function'' introduced in \cite{IkL}.
Therein inverse obstacle scattering problems at a fixed frequency in two dimensions are considered.   
Needless to say, formula (1.17) is not derived in \cite{Ik}.

The condition (1.15) is {\it stable} with respect to the parturbation of $\vartheta\in S(\omega)$ since from the expression (1.3)
we see that the function $S(\omega)\ni\vartheta\longmapsto C_{(p(\omega),\,\omega)}(\delta,Q,\,\vartheta)$ 
is continuous, where the topology of $S(\omega)$ is the relative one from $\Bbb R^3$.  This fact yields a corollary as follows.

\proclaim{\noindent Corollary 1.1.}
Let $\omega$ be regular with respect to $D$.
Under the same assumptions as those in Theorem 1.2 the point $p(\omega)$ is uniquely determined by
the Cauchy data of $u$ on $\partial\Omega$.

\endproclaim

{\it\noindent Proof.}  From (1.16) one has $h_D(\omega)=p(\omega)\cdot\omega$.  
Choose $\vartheta'\in S(\omega)$ sufficiently near $\vartheta$ in such a way that
$C_{(p(\omega),\omega)}(\delta,Q,\vartheta')\not=0$.  Then from the formula (1.17) for two linearly independent directions $\vartheta$ and $\vartheta'$
one gets $p(\omega)\cdot\vartheta$ and $p(\omega)\cdot\vartheta'$.

\noindent
$\Box$

As another direct corollary of Theorem 1.2 and Proposition 1.2 in the case $n=3$ we have the following result.
\proclaim{\noindent Corollary 1.2.}
Assume that $D$ is given by the inside of a convex polyhedron and in a neighbourhood of each vertex $p$ of $D$,  
the $D$ coincides with the inside of a tetrahedron with apex $p$ and that the source $F=F_{\rho, D}$ given by (1.7)
is active at $p$.
Then, we have all the formulae (1.16), (1.17) and (1.18) for all $\omega$ regular with respect to $D$
and $\vartheta\in S(\omega)$.
\endproclaim
{\it\noindent Proof.}  We have: $D$ has a conical singularity from the direction $\omega$ that is regular with respect to $D$
with a triangle $Q$ at each $p(\omega)$.  Thus (1.15) is valid for all $\omega$ regular with respect to $D$ and $\vartheta\in S(\omega)$.
Therefore, we have all the formulae (1.16), (1.17) and (1.18) for all $\omega$ regular with respect to $D$ and $\vartheta\in S(\omega)$.
\noindent
$\Box$

{\bf\noindent Remark 1.1.}  Under the same assumptions as Corollary 1.2 one gets a uniqueness theorem: the Cauchy data of $u$ on $\partial\Omega$ uniquely determines $D$.
The proof is as follows.  From (1.16) one gets $h_D(\omega)$ for all $\omega$ regular with respect to $D$.  The set of all $\omega$ that are not regular with respect to $D$
consists of a set of finite points and arcs on $S^2$.  
This yields the set of all $\omega$ that are regular with respect to $D$ is dense and thus one gets 
$h_D(\omega)$ for all $\omega\in S^2$ because of the continuity of $h_D$.  Therefore one obtains the convex hull of $D$ and thus $D$ itself by the convexity assumption.
This proof is remarkable and unique since 
we never make use of the {\it traditional contradiction argument}`` Suppose we have two different source domains
$D_1$ and $D_2$ which yields the same Cauchy data,...''; any {\it unique continuation argument} of the solution of the governing equation.
One can see such two arguments in \cite{N} in the case when $k=0$ for an inverse problem for detecting a source of {\it gravity anomaly}.

Some of typical examples of $D$ covered by Corollary 1.2 are tetrahedron, regular hexahedron (cube), regular dodecahedron.

So now the central problem in applying Theorem 1.2 to Problem 1 for the source with various source domain 
under our framework is to clarify
the condition (1.15) for general $Q$.
In contrast to Proposition 1.2, when $Q$ is general, we do not know whether there exists a unit vector $\vartheta\in S(\omega)$ such that (1.15) is valid
or not.  Going back to (1.3), we have an explicit vector equation for the constant $C_{(p,\omega)}(\delta,Q,\vartheta)$, if $Q$ is given by the inside of a polygon.
See Proposition 4 in \cite{IkC}.  However, comparing with the case when $Q$ is given by the inside of a triangle,
it seems difficult to deduce the non-vanishing $C_{(p,\omega)}(\delta,Q,\vartheta)$
for all $\vartheta\in S(\omega)$ from the equation directly.  This is an open problem.

\subsection{Explicit formula and its implication}

In this paper, instead of considering general $Q$, we consider another special
$Q$.  It is the case when $Q$ is given by the section of the inside of a {\it circulr cone} by a plane.

Given $p\in\Bbb R^3$, $\mbox{\boldmath $n$}\in S^2$ and $\theta\in\,]0,\,\frac{\pi}{2}[$ let $V_p(-\mbox{\boldmath $n$},\theta)$
denote the inside of the {\it circular cone} with {\it apex} at $p$ and the opening angle $\theta$ around the direction $-\mbox{\boldmath $n$}$,
that is
$$\displaystyle
V_p(-\mbox{\boldmath $n$},\theta)=\left\{x\in\Bbb R^3\,\left\vert\right.
\,(x-p)\cdot(-\mbox{\boldmath $n$})>\cos\theta\,\right\}.
$$
Given $\omega\in S^2$ set
$$\displaystyle
Q=\mbox{\boldmath $V$}_p(-\mbox{\boldmath $n$},\theta)
\cap\left\{x\in\Bbb R^3\,\left\vert\right.\,x\cdot\omega=p\cdot\omega-\delta\,\right\}.
\tag {1.19}
$$
To ensure that $Q$ is non empty and bounded, we impose the restriction between $\omega$ and $\mbox{\boldmath $n$}$ as follows:
$$
\omega\cdot\mbox{\boldmath $n$}>\cos(\pi/2-\theta)=\sin\theta(>0).
\tag{1.20}
$$
This means that the angle between $\omega$ and $\mbox{\boldmath $n$}$ has to be less than $\frac{\pi}{2}-\theta$.
Then it is known that $Q$ is an ellipse and we have
$$\displaystyle
D_{(p,\omega)}(\delta, Q)=\mbox{\boldmath $V$}_p(-\mbox{\boldmath $n$},\theta)
\cap\left\{x\in\Bbb R^3\,\left\vert\right.\,x\cdot\omega>p\cdot\omega-\delta\,\right\}.
\tag {1.21}
$$
The problem here is to compute the complex constant
$C_{(p,\omega)}(\delta,Q,\vartheta)$ with all $\vartheta\in S(\omega)$ for this domain
$D_{(p,\omega)}(\delta,Q)$ with $Q$ given by (1.19).

Instead of (1.3) we employ
the formula (1.4) with $D=D_{(p,\omega)}(\delta,Q)$ with $n=3$:
$$\displaystyle
C_{(p,\omega)}(\delta,Q,\vartheta)
=
\lim_{\tau\longrightarrow\infty}\tau^3e^{-\tau p\cdot(\omega+i\vartheta)}
\int_{D_{(p,\omega)}(\delta,Q)}\,e^{\tau x\cdot(\omega+i\vartheta)}dx.
\tag {1.22}
$$
Here we rewrite this formula.  Choosing sufficiently small positive numbers $\delta'$ and $\delta''$ with $\delta''<\delta'$,
we see that the set 
$$\displaystyle
D_{(p,\omega)}(\delta, Q)\cap\left\{x\in\Bbb R^3\,\left\vert\right.\,x\cdot\mbox{\boldmath $n$}<p\cdot\mbox{\boldmath $n$}-\delta'\,\right\}
$$
is containted in the half-space $x\cdot\omega<p\cdot\omega-\delta''$.

This yields
$$
\displaystyle
e^{-\tau p\cdot(\omega+i\vartheta)}
\int_{D_{(p,\omega)}(\delta,Q)}\,e^{\tau x\cdot(\omega+i\vartheta)}dx
=e^{-\tau p\cdot(\omega+i\vartheta)}
\int_{V}\,e^{\tau x\cdot(\omega+i\vartheta)}dx+O(e^{-\tau\delta''}),
$$
where
$$\displaystyle
V=\mbox{\boldmath $V$}_p(-\mbox{\boldmath $n$},\theta)
\cap\left\{x\in\Bbb R^3\,\left\vert\right.\,x\cdot\mbox{\boldmath $n$}>p\cdot\mbox{\boldmath $n$}-\delta'\,\right\}.
$$
Thus from (1.22) we obtain a more convenient expression
$$\displaystyle
C_{(p,\omega)}(\delta,Q,\vartheta)
=
\lim_{\tau\longrightarrow\infty}\tau^3e^{-\tau p\cdot(\omega+i\vartheta)}
\int_{V}\,e^{\tau x\cdot(\omega+i\vartheta)}dx.
\tag {1.23}
$$
Using this expression we have the following explicit formula of $C_{(p,\omega)}(\delta,Q,\vartheta)$ for $D_{(p,\omega)}(\delta,Q)$
given by (1.21).

\proclaim{\noindent Proposition 1.3.}
We have
$$\displaystyle
C_{(p,\omega)}(\delta, Q,\vartheta)
=6\,V(\theta)\,
(\mbox{\boldmath $n$}\cdot(\omega+i\vartheta)\,)^{-3},
\tag {1.24}
$$
where
$$\displaystyle
V(\theta)=\frac{\pi}{3}\cos\,\theta\sin^2\,\theta.
$$
\endproclaim

Note that the value $V(\theta)$ coincides with 
the volume of the circular cone with the height $\cos\theta$ and the opening angle $\theta$.
This function of $\theta\in\,]0,\,\frac{\pi}{2}\,[$ is monotone increasing in $]0,\,\tan^{-1}\sqrt{2}[$
and decreasing in $]\tan^{-1}\sqrt{2},\,\frac{\pi}{2}[$; takes the maximum value
$\frac{2\pi}{9\sqrt{3}}$ at $\theta=\tan^{-1}\sqrt{2}$.

Now we describe an application to Problem 1.
First we introduce a singularity of a circular cone type for the source domain.

{\bf\noindent Definition 1.3.}  Let $D$ be a non empty bounded open set of $\Bbb R^3$. 
Let $p\in\partial D$.  We say that $D$ has a {\it circular cone singularity} at $p$ if there exist
a positive number $\epsilon$, unit vector $\mbox{\boldmath $n$}$ and number $\theta\in\,]0,\,\frac{\pi}{2}[$ such that
$$\displaystyle
D\cap B_{\epsilon}(p)=V_{p}(-\mbox{\boldmath $n$},\theta)\cap B_{\epsilon}(p).
$$

It is easy to see that notion of the circular cone singularity is a special case of that of the conical one
in the following sense.

\proclaim{\noindent Lemma 1.1.}
Let $\omega\in S^2$ be regular with respect to $D$.
Assume that $D$ has a circular cone singularity at $p(\omega)$.  Then, $D$ has a conical singularity from direction $\omega$
at $p(\omega)$.  More precisely, for a sufficiently small $\delta$ we have the expression
$$\displaystyle
D\cap\left\{x\in\Bbb R^3\,\vert\,  h_D(\omega)-\delta<x\cdot\omega<h_D(\omega)\,\right\}
=D_{(p(\omega),\omega)}(\delta, Q),
$$
where $Q$ is given by (1.19) with $V_{p}(-\mbox{\boldmath $n$},\theta)$ at $p=p(\omega)$ in the definition 1.3 satisfying (1.20).

\endproclaim

As a diect corollary of Theorems 1.1-1.2, Proposition 1.3 and Lemma 1.1 we immediately 
obtain all the results in Theorem 1.2 without the condition (1.15).  We suumarize one of the result as Corollary 1.3 as follows.

\proclaim{\noindent Corollary 1.3(Detecting the point $p(\omega)$).}  Let $u\in H^1(\Omega)$ be an arbitrary solution of (1.6) with the source $F=F_{\rho,D}$ given by (1.7).
Let $\omega\in S^2$ be regular with respect to $D$.  Assume that: $D$ has a circular cone singularity
at $p=p(\omega)$; the source $F$ is active at $p(\omega)$.
Choose two linearly independent vectors $\vartheta=\vartheta_1$ and $\vartheta_2$ in $S(\omega)$.
Then, the point $p(\omega)$ itself and thus $h_D(\omega)=p(\omega)\cdot\omega$ can be extracted from the Cauchy data of $u$ on $\partial\Omega$
by using the formula
$$\displaystyle
p(\omega)\cdot\omega+i\,p(\omega)\cdot\vartheta_j
=\lim_{\tau\rightarrow\infty}
\frac{I_{\omega,\vartheta_j}'(\tau)}{I_{\omega,\vartheta_j}(\tau)},\,\,\,j=1,2.
\tag {1.25}
$$
\endproclaim

By virtue of the formula (1.24), the function $I(\omega,\,\cdot\,)$ has the expression
$$\displaystyle
I(\omega,\vartheta)=6\,\tilde{\rho}(p(\omega))\,V(\theta)(\mbox{\boldmath $n$}\cdot(\omega+i\vartheta))^{-3}.
\tag {1.26}
$$
Formula (1.26) yields the following results.

\proclaim{\noindent Corollary 1.4.}  Let $u\in H^1(\Omega)$ be a solution of (1.6) with the source $F=F_{\rho,D}$ given by (1.7).
Let $\omega\in S^2$ be regular with respect to $D$.  Assume that: $D$ has a circular cone singularity
at $p(\omega)$ such as $D\cap B_{\epsilon}(p(\omega))=V_{p(\omega)}(-\mbox{\boldmath $n$},\theta)\cap B_{\epsilon}(p(\omega))$ with a $\epsilon>0$.

\noindent
(i)  Assume that $F$ is active at $p(\omega)$.    
The vector $\omega$ coincides with $\mbox{\boldmath $n$}$ if and only if the function
$I(\omega,\,\cdot\,)$ is a constant function.

\noindent
(ii)  The vector $\mbox{\boldmath $n$}$ and $\theta$ of $V_{p(\omega)}(-\mbox{\boldmath $n$},\theta)$
and the source strength $\tilde{\rho}(p(\omega))$ satisfies the following two equations:
$$\displaystyle
6\,\vert\tilde{\rho}(p(\omega))\vert\,V(\theta)=(\mbox{\boldmath $n$}\cdot\omega)^3
\max_{\vartheta\in S(\omega)}\vert I(\omega,\vartheta)\vert;
\tag {1.27}
$$
$$\displaystyle
6\,\tilde{\rho}(p(\omega))
\,V(\theta)\,(3(\mbox{\boldmath $n$}\cdot\omega)^2-1)
=\frac{1}{\pi}\,\int_{S(\omega)}\,I(\omega,\vartheta)
\,ds(\vartheta).
\tag {1.28}
$$

\endproclaim

Using the equations (1.26), (1.27) and (1.28) one gets the following corollary.

\proclaim{\noindent Corollary 1.5.}  Let $u\in H^1(\Omega)$ be a solution of (1.6) with the source $F=F_{\rho,D}$ given by (1.7).
Let $\omega\in S^2$ be regular with respect to $D$.  Assume that: $D$ has a circular cone singularity
at $p(\omega)$ such as $D\cap B_{\epsilon}(p(\omega))=V_{p(\omega)}(-\mbox{\boldmath $n$},\theta)\cap B_{\epsilon}(p(\omega))$ with a $\epsilon>0$.
Assume that $F$ is active at $p(\omega)$ and that $\omega\approx\mbox{\boldmath $n$}$ in the sense that
$$\displaystyle
\mbox{\boldmath $n$}\cdot\omega>\frac{1}{\sqrt{3}}.
\tag {1.29}
$$
Then, the value $\gamma=\mbox{\boldmath $n$}\cdot\omega$ is the unique solution of the following quintic equation in $]\,\frac{1}{\sqrt{3}},\,1]$:
$$\displaystyle
\gamma^3(3\gamma^2-1)=
\frac{\displaystyle\left\vert\int_{S(\omega)}\,I(\omega,\vartheta)
\,ds(\vartheta)\right\vert}{\pi\,\max_{\vartheta\in S(\omega)}\vert I(\omega,\vartheta)\vert}.
\tag {1.30}
$$
Besides, for an arbitrary $\vartheta\in S(\omega)$ the value $\mu=\mbox{\boldmath $n$}\cdot\vartheta$ is given by the formulae
$$
\displaystyle
\mu^2=\frac{\displaystyle\gamma^3-\text{Re}\,T(\omega,\vartheta)}
{3\gamma}
\tag {1.31}
$$
and
$$\displaystyle
\mu=\frac{\displaystyle\text{Im}\,T(\omega,\vartheta)}{3\gamma^2-\mu^2},
\tag {1.32}
$$
where
$$\displaystyle
T(\omega,\vartheta)
=\frac{\displaystyle
\int_{S(\omega)}\,I(\omega,\vartheta)\,ds(\vartheta)}
{\pi(3\gamma^2-1)I(\omega,\vartheta)}.
\tag{1.33}
$$

\endproclaim

The condition (1.29) is equivalent to the statement\footnote{We have
$$\displaystyle
\frac{3\pi}{10}+\frac{\pi}{100}>\tan^{-1}\sqrt{2}>\frac{3\pi}{10}.
$$
}: the angle between $\omega$ and $\mbox{\boldmath $n$}$
is less than $\tan^{-1}\sqrt{2}$.  Thus it is not so strict condition.  
The denominator of (1.32) is not zero because of $3\gamma^2-\mu^2\ge 3\gamma^2-1$ and (1.29).

Under the same assumptions as Corollary 1.5, one can finally calculate the quantity
$$\displaystyle
\tilde{\rho}(p(\omega))\,V(\theta)
\tag {1.34}
$$
and $\mbox{\boldmath $n$}$ from the Cauchy data of $u$ on $\partial\Omega$.  This is the final conclusion.

The procedure is as follows.

\noindent
{\bf Step 1.}  Calculate $p(\omega)$ via the formula (1.25).

\noindent
{\bf Step 2.}  Calculate $I(\omega,\vartheta)$ via the formula (1.14) and the computed $p(\omega)$ in Step 1.

\noindent
{\bf Step 3.}  If $I(\omega,\vartheta)$ looks like a constant function, decide $\omega\approx\mbox{\boldmath $n$}$ in the sense (1.29).
If not so, search another $\omega$ around the original one in such a way that $\omega\approx\mbox{\boldmath $n$}$ as above by try and error
and finally fix it.

\noindent
{\bf Step 4.}  Find the value $\gamma=\mbox{\boldmath $n$}\cdot\omega$ by solving the quintic equation (1.30).

\noindent
{\bf Step 5.}  Find the value (1.34) via the formulae (1.28) with the computed $\mbox{\boldmath $n$}\cdot\omega$
in Step 4.

\noindent
{\bf Step 6.}  Choose linearly independent vectors $\vartheta_1, \vartheta_2\in S(\omega)$ and calculate
$T(\omega,\vartheta_j)$, $j=1,2$ via the formula (1.33) using the computed value $\gamma$ in Step 4.

\noindent
{\bf Step 7.}  Find $\mu=\mu_j=\mbox{\boldmath $n$}\cdot\vartheta_j$ by solving (1.31) and (1.32) using the computed $T(\omega,\vartheta_j)$
in Step 6.

\noindent
{\bf Step 8.}  Find $\mbox{\boldmath $n$}$ by solving $\mbox{\boldmath $n$}\cdot\omega=\gamma$, $\mbox{\boldmath $n$}\cdot\vartheta_j=\mu_j$, $j=1,2$.

Note that, in addition, if the opening angle $\theta$/the source strength $\tilde{\rho}(p(\omega))$ is known, then one obtains the value of
$\tilde{\rho}(p(\omega))$/the volume $V(\theta)$ via the 
computed value (1.34) in Step 5.

This paper is organized as follows.  In the next section we give a proof of Proposition 1.3.  
It is based on the integral representation (2.8) of the complex constant $C_{(p,\omega)}(\delta, Q,\vartheta)$ and the residue calculus.
Proofs of Corollaries 1.4 and 1.5 are given in Section 3.
In Section 4, an inverse obstacle problem for a penetrable obstacle in three dimensions is considered.
The corresponding results in this case are given and in Section 5 a possible direction of the extension
of all the results in this paper is commented.  Appendix is devoted to an example covered by the results in Section 4.

\section{Proof of Proposition 1.3}

In order to compute this right-hand side, we choose two unit vectors $\mbox{\boldmath $l$}$ and $\mbox{\boldmath $m$}$ perpendicular
to each other in such a way that $\mbox{\boldmath $n$}=\mbox{\boldmath $l$}\times\mbox{\boldmath $m$}$.

We see that the intersection of $\partial V_p(-\mbox{\boldmath $n$},\theta)$ with the plane
$(x-p)\cdot\mbox{\boldmath $n$}=-(1/\tan\,\theta)$ coincides with the circle with radius $1$ centered at the point
$p-(1/\tan\,\theta)\mbox{\boldmath $n$}$ on the plane.
The pointing vector of an arbitrary point on the circle with respect to point $p$ has the expression
$$\displaystyle
\vartheta(w)=\cos\,w\,\mbox{\boldmath $l$}+\sin\,w\,\mbox{\boldmath $m$}
-\frac{1}{\tan\,\theta}\,\mbox{\boldmath $n$}
\tag {2.1}
$$
with a parameter $w\in\,[0,2\pi]$.  
Besides, from the geometrical meaning of $\vartheta(w)$, we have 
$$
\displaystyle\max_{w\in[0,\,2\pi]}\,\vartheta(w)\cdot\omega<0.
\tag {2.2}
$$

\proclaim{\noindent Lemma 2.1.}
We have the expression
$$\displaystyle
(\omega+i\vartheta)\,C_{(p,\omega)}(\delta,Q,\vartheta)
=\frac{1}{\tan\,\theta}
\int_0^{2\pi}
\frac{\cos\,w\,\mbox{\boldmath $l$}+\sin\,w\,\mbox{\boldmath $m$}+\tan\,\theta\,\mbox{\boldmath $n$}}
{\{\vartheta (w)\cdot(\omega+i\vartheta)\}^2}dw.
\tag {2.3}
$$

\endproclaim

{\it\noindent Proof.}
Let $\mbox{\boldmath $a$}$ be an arbitrally three dimensional complex vector.
We have
$$\displaystyle
\int_{V}\,\nabla\cdot(e^{\tau x\cdot(\omega+i\vartheta)}\mbox{\boldmath $a$})\,dx
=\tau(\omega+i\vartheta)\cdot\mbox{\boldmath $a$}\,\int_{V}\,e^{\tau x\cdot(\omega+i\vartheta)}dx.
$$
The divergence theorem yields
$$\displaystyle
(\omega+i\vartheta)\cdot\mbox{\boldmath $a$}\,\int_{V}\,e^{\tau x\cdot(\omega+i\vartheta)}dx
=\tau^{-1}\int_{\partial V}\,e^{\tau x\cdot(\omega+i\vartheta)}\mbox{\boldmath $a$}\cdot\mbox{\boldmath $\nu$}\,dS(x),
\tag {2.4}
$$
where $\mbox{\boldmath $\nu$}$ denotes the outer unit normal vector to $\partial V$.

Decompose $\partial V=V_1\cup V_2$ with $V_1\cap V_2=\emptyset$, where
$$\begin{array}{l}
\displaystyle
V_1=\{x\,\vert\,-(x-p)\cdot\mbox{\boldmath $n$}=\vert x-p\vert\cos\,\theta,\,
-\delta'<(x-p)\cdot\mbox{\boldmath $n$}<0\},\\
\\
\displaystyle
V_2=\{x\,\vert\,\vert x-(p-\delta'\,\mbox{\boldmath $n$})\vert\le\delta'\,\tan\,\theta,\,
(x-p)\cdot\mbox{\boldmath $n$}=-\delta'\}.
\end{array}
$$

To compute the surface integral over $V_1$, we make use of the change of variables as follows:
$$\begin{array}{ll}
\displaystyle
x
&
\displaystyle
=(p-\delta'\,\mbox{\boldmath $n$})+r(\cos\,w\,\mbox{\boldmath $l$}
+\sin\,w\,\mbox{\boldmath $m$})+\left(\delta'-\frac{r}{\tan\,\theta}\right)\,\mbox{\boldmath $n$}
\\
\\
\displaystyle
&
\displaystyle
=p+r\vartheta(w),
\end{array}
\tag {2.5}
$$
where $(r,w)\in\,[0,\,\delta'\tan\,\theta]\times[0,\,2\pi[$ and $\vartheta(w)$ is given by (2.1).
Then the surface element has the expression
$$\displaystyle
dS(x)=\frac{r}{\sin\,\theta}\,drdw
$$
and outer unit normal $\mbox{\boldmath $\nu$}$ to $V_1$ takes the form
$$\displaystyle
\nu=\sin\,\theta\left(\mbox{\boldmath $n$}+\frac{\cos\,w\,\mbox{\boldmath $l$}+\sin\,w\,\mbox{\boldmath $m$}}{\tan\,\theta}
\right).
$$

Now from (2.4) and the decomposition $\partial V=V_1\cup V_2$, we have
$$\begin{array}{l}
\displaystyle
\,\,\,\,\,\,
e^{-\tau p\cdot(\omega+i\vartheta)}
(\omega+i\vartheta)\cdot\mbox{\boldmath $a$}\int_{V} v\,dx\\
\\
\displaystyle
=e^{-\tau p\cdot(\omega+i\vartheta)}\tau^{-1}
\int_{V_1} v\mbox{\boldmath $a$}\cdot\nu\,dS(x)
-e^{-\tau p\cdot(\omega+i\vartheta)}\tau^{-1}
\int_{V_2} v\mbox{\boldmath $a$}\cdot\mbox{\boldmath $n$}\,dS(x)\\
\\
\displaystyle
\equiv I+II,
\end{array}
\tag {2.6}
$$
where $v=e^{\tau x\cdot(\omega+i\vartheta)}$.

Since the set $V_2$ is contained in the half-space $x\cdot\omega\le p\cdot\omega-\delta''$, one gets
$$
\displaystyle
II=O(\tau^{-1}e^{-\tau\delta''}).
\tag {2.7}
$$

On $I$, using the change of variables given by (2.5), one has
$$\begin{array}{c}
\displaystyle
x\cdot\omega=p\cdot\omega+r\,\vartheta(w)\cdot\omega,\\
\\
\displaystyle
x\cdot\vartheta=p\cdot\vartheta+r\,\vartheta(w)\cdot\vartheta.
\end{array}
$$
And also noting (2.2), one gets
$$\begin{array}{ll}
\displaystyle
\tau\,I
&
\displaystyle
=\int_0^{2\pi}dw
\int_0^{\delta'\tan\,\theta}
rdr
e^{\tau r\vartheta(w)\cdot\omega+i\tau\,r\vartheta(w)\cdot\vartheta}
\left(\mbox{\boldmath $n$}+\frac{\cos\,w\,\mbox{\boldmath $l$}
+\sin\,w\,\mbox{\boldmath $m$}}{\tan\,\theta}\right)\cdot\mbox{\boldmath $a$}\\
\\
\displaystyle
&
\displaystyle
=\frac{1}{\tau^2}
\int_0^{2\pi}dw
\int_0^{\tau\delta'\tan\,\theta}
sds
e^{s\vartheta(w)\cdot\omega+i\,s\vartheta(w)\cdot\vartheta}
\left(\mbox{\boldmath $n$}+\frac{\cos\,w\,\mbox{\boldmath $l$}
+\sin\,w\,\mbox{\boldmath $m$}}{\tan\,\theta}\right)\cdot\mbox{\boldmath $a$}\\
\\
\displaystyle
&
\displaystyle
=\frac{1}{\tau^2}
\int_0^{2\pi}dw
\int_0^{\infty}
sds
e^{s\vartheta(w)\cdot\omega+is\vartheta(w)\cdot\vartheta}
\left(\mbox{\boldmath $n$}+\frac{\cos\,w\,\mbox{\boldmath $l$}
+\sin\,w\,\mbox{\boldmath $m$}}{\tan\,\theta}\right)\cdot\mbox{\boldmath $a$}+O(\tau^{-4}).
\end{array}
$$
Here one can 
apply the following formula to this right-hand side:
$$\displaystyle
\int_0^{\infty}se^{as}e^{ibs}ds=\frac{1}{(a+ib)^2},\,\,a<0.
$$
Then one gets
$$\displaystyle
I=\frac{1}{\tau^3\tan\,\theta}
\int_0^{2\pi}
\frac{\cos\,w\,\mbox{\boldmath $l$}
+\sin\,w\,\mbox{\boldmath $m$}+\tan\,\theta\,\mbox{\boldmath $n$}}
{\{\vartheta(w)\cdot(\omega+i\Theta)\}^2}\,dw+O(\tau^{-5}).
$$
Now this together with (1.23), (2.6) and (2.7) yields the desired formula.

\noindent
$\Box$

Now from (1.20) and (2.3) we have the integral representation of $C_{(p,\omega)}(\delta,Q,\vartheta)$:
$$\displaystyle
C_{(p,\omega)}(\delta,Q,\vartheta)
=\frac{1}{\mbox{\boldmath $n$}\cdot(\omega+i\vartheta)}
\int_0^{2\pi}\frac{dw}{\{\vartheta(w)\cdot(\omega+i\vartheta)\}^2}.
\tag {2.8}
$$
This formula shows that the constant $C_{(p,\omega)}(\delta,Q,\vartheta)$ is independent
of $p$ and $\delta$ when $Q$ is given by (1.19).

By computing the integral of the right-hand side on (2.8) we obtain the explicit value of $C_{(p,\omega}(\delta, Q,\vartheta)$.

\proclaim{\noindent Lemma 2.2.}  We have: $C_{(p,\omega)}(\delta, Q,\vartheta)\not=0$ if and only if
$$\displaystyle
\frac{\sin\,\theta}
{1+\cos\,\theta}<\left\vert
\frac{\mbox{\boldmath $n$}\cdot(\omega+i\vartheta)}
{(\mbox{\boldmath $l$}-i\mbox{\boldmath $m$})\cdot(\omega+i\vartheta)}
\right\vert<\frac{1+\cos\,\theta}{\sin\,\theta}
\tag {2.9}
$$
and then
$$\displaystyle
C_{(p,\omega)}(\delta, Q,\vartheta)
=2\pi\cos\theta\,\sin^2\,\theta\,
(\mbox{\boldmath $n$}\cdot(\omega+i\vartheta)\,)^{-3}.
\tag {2.10}
$$

\endproclaim

{\it\noindent Proof.}
Set
$$\displaystyle
A=\mbox{\boldmath $l$}\cdot(\omega+i\vartheta), \,\,B=\mbox{\boldmath $m$}\cdot(\omega+i\vartheta),\,\,
C=-\frac{1}{\tan\,\theta}\mbox{\boldmath $n$}\cdot(\omega+i\vartheta)
$$
and $z=e^{iw}$.  One can write
$$\begin{array}{ll}
\displaystyle
\vartheta(w)\cdot(\omega+i\vartheta)
&
\displaystyle
=A\cos\,w+B\sin\,w+C\\
\\
\displaystyle
&
\displaystyle
=\frac{A}{2}(z+z^{-1})-i\frac{B}{2}(z-z^{-1})+C\\
\\
\displaystyle
&
\displaystyle
=\frac{1}{2z}
\{(A-iB)z^2+2Cz+(A+iB)\}.
\end{array}
$$
  
Here we claim 
$$
\displaystyle
A-iB\equiv(\mbox{\boldmath $l$}-i\mbox{\boldmath $m$})\cdot(\omega+i\vartheta)\not=0.
\tag {2.11}
$$

Assume contrary that $A-iB=0$.
Since we have $$\displaystyle
A-iB
=\mbox{\boldmath $l$}\cdot\omega+\mbox{\boldmath $m$}\cdot\vartheta
+i(\mbox{\boldmath $l$}\cdot\vartheta-\mbox{\boldmath $m$}\cdot\omega),
$$
it must hold that
$$\displaystyle
\mbox{\boldmath $l$}\cdot\omega=-\mbox{\boldmath $m$}\cdot\vartheta,\,\,
\mbox{\boldmath $m$}\cdot\omega=\mbox{\boldmath $l$}\cdot\vartheta.
\tag {2.12}
$$
Then we have
$$\begin{array}{ll}
\displaystyle
(\mbox{\boldmath $n$}\cdot\vartheta)^2
&
\displaystyle
=\vert\vartheta\vert^2-(\mbox{\boldmath $l$}\cdot\vartheta)^2-(\mbox{\boldmath $m$}\cdot\vartheta)^2
\\
\\
\displaystyle
&
\displaystyle
=\vert\omega\vert^2-(\mbox{\boldmath $l$}\cdot\omega)^2-(\mbox{\boldmath $m$}\cdot\omega)^2
\\
\\
\displaystyle
&
\displaystyle
=(\mbox{\boldmath $n$}\cdot\omega)^2.
\end{array}
\tag {2.13}
$$

On the other hand, we have
$$\displaystyle
0=\omega\cdot\vartheta
=(\mbox{\boldmath $l$}\cdot\omega)(\mbox{\boldmath $l$}\cdot\vartheta)
+(\mbox{\boldmath $m$}\cdot\omega)(\mbox{\boldmath $m$}\cdot\vartheta)
+(\mbox{\boldmath $n$}\cdot\omega)(\mbox{\boldmath $n$}\cdot\vartheta).
$$
Here by (2.12) one has
$(\mbox{\boldmath $l$}\cdot\omega)(\mbox{\boldmath $l$}\cdot\vartheta)
+(\mbox{\boldmath $m$}\cdot\omega)(\mbox{\boldmath $m$}\cdot\vartheta)=0$.
Thus one obtains 
$$
\displaystyle
0=(\mbox{\boldmath $n$}\cdot\omega)(\mbox{\boldmath $n$}\cdot\vartheta).
$$
Now a combination of this and (2.13) yields $\mbox{\boldmath $n$}\cdot\omega=0$.
However, by (1.20) this is impossible.

Therefore we obtain the expression
$$\displaystyle
\vartheta(w)\cdot(\omega+i\vartheta)
=\frac{A-iB}{2z}f(z)\vert_{z=e^{iw}},
\tag {2.14}
$$
where
$$\displaystyle
f(z)=
\left(z+\frac{C}{A-iB}\right)^2
-\frac{C^2-(A^2+B^2)}{(A-iB)^2}.
$$
Here we write
$$\begin{array}{ll}
\displaystyle
C^2-(A^2+B^2)
&
\displaystyle
=\frac{1}{\tan^2\,\theta}
(\mbox{\boldmath $n$}\cdot(\omega+i\vartheta))^2
-\{(\mbox{\boldmath $l$}\cdot(\omega+i\vartheta))^2
+(\mbox{\boldmath $m$}\cdot(\omega+i\vartheta))^2\}\\
\\
\displaystyle
&
\displaystyle
=
\frac{1}{\tan^2\,\theta}
\{(\mbox{\boldmath $n$}\cdot\omega)^2-(\mbox{\boldmath $n$}\cdot\vartheta)^2
+2i(\mbox{\boldmath $n$}\cdot\omega)(\mbox{\boldmath $n$}\cdot\vartheta\}\\
\\
\displaystyle
&
\displaystyle
\,\,\,
-\{(\mbox{\boldmath $l$}\cdot\omega)^2
+(\mbox{\boldmath $m$}\cdot\omega)^2
-(\mbox{\boldmath $l$}\cdot\vartheta)^2
-(\mbox{\boldmath $m$}\cdot\vartheta)^2
+2i(\mbox{\boldmath $l$}\cdot\omega)(\mbox{\boldmath $l$}\cdot\vartheta)
+2i(\mbox{\boldmath $m$}\cdot\omega)(\mbox{\boldmath $m$}\cdot\vartheta)\}\\
\\
\displaystyle
&
\displaystyle
=\left(\frac{1}{\tan^2\,\theta}+1\right)(\mbox{\boldmath $n$}\cdot\omega)^2
-
\left(\frac{1}{\tan^2\,\theta}+1\right)(\mbox{\boldmath $n$}\cdot\vartheta)^2
+\frac{1}{\tan^2\,\theta}\,2i(\mbox{\boldmath $n$}\cdot\omega)(\mbox{\boldmath $n$}\cdot\vartheta)\\
\\
\displaystyle
&
\displaystyle
\,\,\,
-2i\{
(\mbox{\boldmath $l$}\cdot\omega)(\mbox{\boldmath $l$}\cdot\vartheta)
+(\mbox{\boldmath $m$}\cdot\omega)(\mbox{\boldmath $m$}\cdot\vartheta)\}\\
\\
\displaystyle
&
\displaystyle
=\frac{1}{\sin^2\,\theta}\{(\mbox{\boldmath $n$}\cdot\omega)^2
-(\mbox{\boldmath $n$}\cdot\vartheta)^2\}
+\frac{1}{\sin^2\,\theta}\,2i(\mbox{\boldmath $n$}\cdot\omega)(\mbox{\boldmath $n$}\cdot\vartheta)
\\
\\
\displaystyle
&
\displaystyle
\,\,\,
-2i\{
(\mbox{\boldmath $l$}\cdot\omega)(\mbox{\boldmath $l$}\cdot\vartheta)
+(\mbox{\boldmath $m$}\cdot\omega)(\mbox{\boldmath $m$}\cdot\vartheta)
+(\mbox{\boldmath $n$}\cdot\omega)(\mbox{\boldmath $n$}\cdot\vartheta)\}
\\
\\
\displaystyle
&
\displaystyle
=\frac{1}{\sin^2\,\theta}\{(\mbox{\boldmath $n$}\cdot\omega)^2
-(\mbox{\boldmath $n$}\cdot\vartheta)^2\}
+\frac{1}{\sin^2\,\theta}\,2i(\mbox{\boldmath $n$}\cdot\omega)(\mbox{\boldmath $n$}\cdot\vartheta)
\\
\\
\displaystyle
&
\displaystyle
\,\,\,
-2i\omega\cdot\vartheta.
\end{array}
$$
Since $\omega\cdot\vartheta=0$, we finally obtain
$$\displaystyle
C^2-(A^2+B^2)
=\left(\frac{\mbox{\boldmath $n$}\cdot(\omega+i\vartheta)}
{\sin\,\theta}\right)^2.
$$
Now set
$$\displaystyle
z_{\pm}=\frac{(\cos\,\theta\pm 1)}{\sin\,\theta}
\frac{\mbox{\boldmath $n$}\cdot(\omega+i\vartheta)}
{(\mbox{\boldmath $l$}-i\mbox{\boldmath $m$})\cdot(\omega+i\vartheta)}.
\tag {2.15}
$$
Then one gets the factorization
$$\displaystyle
f(z)=(z-z_+)(z-z_{-}).
$$
By (2.15) we have $\vert z_+\vert>\vert z_{-}\vert$.  Besides,  from (2.2), (2.11) and (2.14) we have $f(e^{iw})\not=0$ for all $w\in\,[0,\,2\pi]$.
This ensures that the complex numbers $z_{+}$ and $z_{-}$ are not on the circle $\vert z\vert=1$.

Thus from (2.14) one gets
$$
\displaystyle
\int_0^{2\pi}
\frac{dw}
{\{\vartheta(w)\cdot(\omega+i\vartheta)\}^2}
=\frac{4}{i(A-iB)^2}
\int_{\vert z\vert=1}\frac{zdz}{(z-z_{+})^2(z-z_{-})^2}.
\tag {2.16}
$$
The residue calculus yields
$$\displaystyle
\int_{\vert z\vert=1}\frac{zdz}{(z-z_{+})^2(z-z_{-})^2}
=
\left\{
\begin{array}{ll}
\displaystyle
0 & \text{if $\vert z_{-}\vert>1$,}
\\
\\
\displaystyle
0 & \text{if $\vert z_{-}\vert<1$ and $\vert z_{+}\vert<1$,}
\\
\\
\displaystyle
2\pi i\frac{z_{+}+z_{-}}{(z_{+}-z_{-})^3}\not=0
&
\text{if $\vert z_{-}\vert<1<\vert z_{+}\vert$.}
\end{array}
\right.
$$
And also (2.15) gives
$$\begin{array}{ll}
\displaystyle
2\pi i\frac{z_{+}+z_{-}}{(z_{+}-z_{-})^3}
&
\displaystyle
=2\pi i\cdot2\frac{\cos\theta}{\sin\theta}
\frac{\mbox{\boldmath $n$}\cdot(\omega+i\vartheta)}
{(\mbox{\boldmath $l$}-i\mbox{\boldmath $m$})\cdot(\omega+i\vartheta)}
\cdot
(\frac{\sin\theta}{2})^3
\left\{\frac{(\mbox{\boldmath $l$}-i\mbox{\boldmath $m$})\cdot(\omega+i\vartheta)}{\mbox{\boldmath $n$}\cdot(\omega+i\vartheta)}
\right\}^3
\\
\\
\displaystyle
&
\displaystyle
=\frac{\pi i}{2}\cos\,\theta\sin^2\,\theta
\left\{\frac{(\mbox{\boldmath $l$}-i\mbox{\boldmath $m$})\cdot(\omega+i\vartheta)}{\mbox{\boldmath $n$}\cdot(\omega+i\vartheta)}\right\}^2
\\
\\
\displaystyle
&
\displaystyle
=\frac{\pi i}{2}\cos\,\theta\sin^2\,\theta
\left\{\frac{A-iB}
{\mbox{\boldmath $n$}\cdot(\omega+i\vartheta)}\right\}^2.
\end{array}
$$
Thus (2.16) yields
$$
\displaystyle
\int_0^{2\pi}
\frac{dw}
{\{\vartheta(w)\cdot(\omega+i\vartheta)\}^2}
=2\pi\cos\,\theta\sin^2\,\theta
\left\{\frac{1}{\mbox{\boldmath $n$}\cdot(\omega+i\vartheta)}\right\}^2
$$
provided $\vert z_{-}\vert<1<\vert z_{+}\vert$.

From these together with (2.8) we obtain the desired conclusion.

\noindent
$\Box$

Note that (2.10) is nothing but (1.24).
Since (2.9) looks like a condition depending on the choice of $\mbox{\boldmath $l$}$ and $\mbox{\boldmath $m$}$ we further rewrite the number
$$\displaystyle
K(\vartheta;\omega,\mbox{\boldmath $n$})
=\left\vert
\frac{\mbox{\boldmath $n$}\cdot(\omega+i\vartheta)}
{(\mbox{\boldmath $l$}-i\mbox{\boldmath $m$})\cdot(\omega+i\vartheta)}
\right\vert.
$$
We have 
$$\begin{array}{l}
\displaystyle
\,\,\,\,\,\,
\vert(\mbox{\boldmath $l$}-i\mbox{\boldmath $m$})\cdot(\omega+i\vartheta)\vert^2
\\
\\
\displaystyle
=(\mbox{\boldmath $l$}\cdot\omega+\mbox{\boldmath $m$}\cdot\vartheta)^2+
(\mbox{\boldmath $l$}\cdot\vartheta-\mbox{\boldmath $m$}\cdot\omega)^2
\\
\\
\displaystyle
=2-(\mbox{\boldmath $n$}\cdot\omega)^2-(\mbox{\boldmath $n$}\cdot\vartheta)^2
+2(\mbox{\boldmath $l$}\cdot\omega\,\mbox{\boldmath $m$}\cdot\vartheta-\mbox{\boldmath $l$}\cdot\vartheta\,\mbox{\boldmath $m$}\cdot\omega).
\end{array}
$$
Here we see that
$$\displaystyle
\mbox{\boldmath $n$}\cdot(\omega\times\vartheta)
=\mbox{\boldmath $l$}\cdot\omega\,\mbox{\boldmath $m$}\cdot\vartheta-\mbox{\boldmath $l$}\cdot\vartheta\,\mbox{\boldmath $m$}\cdot\omega.
$$
Thus one has
$$\begin{array}{l}
\displaystyle
\,\,\,\,\,\,
\vert(\mbox{\boldmath $l$}-i\mbox{\boldmath $m$})\cdot(\omega+i\vartheta)\vert^2
\\
\\
\displaystyle
=2-(\mbox{\boldmath $n$}\cdot\omega)^2-(\mbox{\boldmath $n$}\cdot\vartheta)^2
+2\mbox{\boldmath $n$}\cdot(\omega\times\vartheta).
\end{array}
$$
Therefore we obtain
$$\displaystyle
K(\vartheta;\omega,\mbox{\boldmath $n$})=\frac{\displaystyle
\sqrt{(\mbox{\boldmath $n$}\cdot\omega)^2+(\mbox{\boldmath $n$}\cdot\vartheta)^2
}}
{\displaystyle
\sqrt{
2-(\mbox{\boldmath $n$}\cdot\omega)^2-(\mbox{\boldmath $n$}\cdot\vartheta)^2
+2\mbox{\boldmath $n$}\cdot(\omega\times\vartheta)}
}.
$$
Besides, we have
$$\displaystyle
\frac{1-\cos\,\theta}{\sin\,\theta}=\tan\,\frac{\theta}{2}
$$
and
$$\displaystyle
\frac{1+\cos\,\theta}{\sin\,\theta}=\frac{1}{\tan\,\frac{\theta}{2}}.
$$
Thus (2.9) is equivalent to the condition
$$\displaystyle
\tan\,\frac{\theta}{2}<
K(\vartheta;\omega,\mbox{\boldmath $n$})
<\frac{1}{\tan\,\frac{\theta}{2}}.
\tag {2.17}
$$

Here consider the case $\mbox{\boldmath $\omega$}\times\mbox{\boldmath $n$}\not=\mbox{\boldmath $0$}$.
Choose
$$\displaystyle
\vartheta=\frac{\mbox{\boldmath $\omega$}\times\mbox{\boldmath $n$}}
{\vert\mbox{\boldmath $\omega$}\times\mbox{\boldmath $n$}\vert}.
$$
We have $\vartheta\cdot\omega=\vartheta\cdot\mbox{\boldmath $n$}=0$ and $\vartheta\in S^2$.

Since we have 
$$\displaystyle
\mbox{\boldmath $n$}\cdot(\omega\times\vartheta)=-\vert\omega\times\mbox{\boldmath $n$}\vert
$$
and
$$\displaystyle
1=(\mbox{\boldmath $n$}\cdot\omega)^2+\vert\omega\times\mbox{\boldmath $n$}\vert^2,
$$
one gets
$$\begin{array}{l}
\displaystyle
\,\,\,\,\,\,
2-(\mbox{\boldmath $n$}\cdot\omega)^2-(\mbox{\boldmath $n$}\cdot\vartheta)^2
+2\mbox{\boldmath $n$}\cdot(\omega\times\vartheta)
\\
\\
\displaystyle
=1+\vert\omega\times\mbox{\boldmath $n$}\vert^2-2\vert\omega\times\mbox{\boldmath $n$}\vert\\
\\
\displaystyle
=(1-\vert\omega\times\mbox{\boldmath $n$}\vert)^2.
\end{array}
$$
Therefore, we obtain
$$\displaystyle
K(\omega\times\mbox{\boldmath $n$};\omega,\mbox{\boldmath $n$})
=\frac{\omega\cdot\mbox{\boldmath $n$}}{1-\vert\omega\times\mbox{\boldmath $n$}\vert}.
$$
Note that we are considering $\omega$ satisfying (1.20).  Let $\varphi$ denote the angle between $\omega$ and $\mbox{\boldmath $n$}$.
Under the condition $\omega\times\mbox{\boldmath $n$}\not=\mbox{\boldmath $0$}$,  we see that (1.20) is equivalent to the condition
$$\displaystyle
0<\varphi<\frac{\pi}{2}-\theta.
\tag {2.18}
$$
Then one can write
$$\begin{array}{ll}
\displaystyle
K(\omega\times\mbox{\boldmath $n$};\omega,\mbox{\boldmath $n$})
&
\displaystyle
=\frac{\cos\,\varphi}{1-\sin\,\varphi}
\\
\\
\displaystyle
&
\displaystyle
=\frac{1+\sin\,\varphi}{\cos\,\varphi}
\\
\\
\displaystyle
&
\displaystyle
=\frac{\displaystyle 1+\cos\,(\frac{\pi}{2}-\varphi)}{\displaystyle
\sin\,(\frac{\pi}{2}-\varphi)}
\\
\\
\displaystyle
&
\displaystyle
=\frac{1}{\displaystyle\tan\,\frac{1}{2}\,(\frac{\pi}{2}-\varphi)}
\end{array}
$$
Thus (2.18) gives
$$\displaystyle
1<K(\omega\times\mbox{\boldmath $n$};\omega,\mbox{\boldmath $n$})<\frac{1}{\displaystyle \tan\,\frac{\theta}{2}}.
\tag {2.19}
$$
Since we have $\tan\,\frac{\theta}{2}<1$ for all $\theta\in\,]0,\,\frac{\pi}{2}[$, (2.19) yields the validity of (2.17).

Next consider the case $\omega\times\mbox{\boldmath $n$}=\mbox{\boldmath $0$}$.  By (1.20) we have $\omega=\mbox{\boldmath $n$}$.
Then, for all $\vartheta$ perpendicular to $\mbox{\boldmath $n$}$ satisfies
$$\displaystyle
K(\vartheta;\mbox{\boldmath $n$}, \mbox{\boldmath $n$})
=1.
$$
This yields that (2.17) is valid for all $\theta\in\,]0,\,\frac{\pi}{2}[$.

The results above are summarized as follows.
Given $\omega\in S^2$ with (1.20) define the subset of $S^2$ 
$$\displaystyle
{\cal K}(\omega;\mbox{\boldmath $n$},\theta)
=
\left\{
\vartheta\in S^2\,\left\vert\right.\,\vartheta\cdot\omega=0,
\,\,\text{$K(\vartheta;\omega,\mbox{\boldmath $n$})$ satisfies (2.19)\,}\,\right\}.
$$
Then, we have

$\bullet$  If $\omega\not\not=\mbox{\boldmath $n$}$, then $\omega\times\mbox{\boldmath $n$}\in{\cal K}(\omega;\mbox{\boldmath $n$},\theta)$.

$\bullet$  If $\omega=\mbox{\boldmath $n$}$, then ${\cal K}(\omega;\mbox{\boldmath $n$},\theta)=
\{\vartheta\in S^2\,\vert\,\vartheta\cdot\omega=0\}\equiv S(\omega)$.

Thus, any way the set ${\cal K}(\omega;\mbox{\boldmath $n$},\theta)$ is not empty and clearly open
with respect to the topology of the set $S(\omega)$ which is the relative topology
of $S^2$.
Besides, we can say more about ${\cal K}(\omega;\mbox{\boldmath $n$},\theta)$.
We claim set ${\cal K}(\omega;\mbox{\boldmath $n$},\theta)$ is closed.
For this, It suffices to show that if a sequence $\{\vartheta_n\}$ of ${\cal K}(\omega;\mbox{\boldmath $n$},\theta)$
converges to a point $\vartheta\in S(\omega)$, then $\vartheta\in{\cal K}(\omega;\mbox{\boldmath $n$},\theta)$.
This is proved as follows.  By assumption, each $\vartheta_n$ satisfies
$$\displaystyle
\tan\,\frac{\theta}{2}<
K(\vartheta_n;\omega,\mbox{\boldmath $n$})
<\frac{1}{\tan\,\frac{\theta}{2}}.
$$
Taking the limit, we have
$$\displaystyle
\tan\,\frac{\theta}{2}\le
K(\vartheta;\omega,\mbox{\boldmath $n$})
\le\frac{1}{\tan\,\frac{\theta}{2}}.
$$
By (2.15) this is equivalent to $\vert z_{+}\vert\ge 1$ and $\vert z_{-}\vert\le 1$.  However, in the proof of
Lemma 2.2 we know that
$\vert z_{+}\vert\not=1$ and $\vert z_{-}\vert\not=1$.  Thus we have $\vert z_{+}\vert>1$ and $\vert z_{-}\vert<1$.
This is equivalent to $\vartheta\in{\cal K}(\omega;\mbox{\boldmath $n$},\theta)$.

Since $S(\omega)$ is connected, ${\cal K}(\omega;\mbox{\boldmath $n$},\theta)$ is not empty, open and closed
we conclude ${\cal K}(\omega;\mbox{\boldmath $n$},\theta)=S(\omega)$.

This completes the proof of Proposition 1.3.

\section{Proof of Corollaries 1.4 and 1.5}

Note that $\omega$ satisfies (1.20).

\subsection{On Corollary 1.4}

From (1.26) we have, if $\omega=\mbox{\boldmath $n$}$, then for all $\vartheta\in S(\omega)$
$$\displaystyle
I(\omega,\vartheta)=
6\tilde{\rho}(p(\omega))\,V(\theta)(\mbox{\boldmath $n$}\cdot\omega)^{-3}.
$$
On the other hand, if $\omega\not=\mbox{\boldmath $n$}$, then we have $\omega\times\mbox{\boldmath $n$}\not=\mbox{\boldmath $0$}$ (under the condition (1.19))
and
$$\displaystyle
S(\omega)\cap S(\mbox{\boldmath $n$})=\left\{\pm\frac{\omega\times\mbox{\boldmath $n$}}{\vert \omega\times\mbox{\boldmath $n$}\vert}\right\}.
$$
Thus one gets
$$\displaystyle
I(\omega,\vartheta)
=
\begin{array}{ll}
\displaystyle
6\tilde{\rho}(p(\omega))\,V(\theta)
\left(\mbox{\boldmath $n$}\cdot\omega\mp\,i\frac{\vert\omega\times\mbox{\boldmath $n$}\vert^2}{\vert\omega\times(\omega\times\mbox{\boldmath $n$})\vert}\right)^{-3}
& \text{for $\displaystyle\vartheta=\pm\frac{\omega\times(\omega\times\mbox{\boldmath $n$})}{\vert\omega\times(\omega\times\mbox{\boldmath $n$})\vert}$.}
\end{array}
$$

Thus one gets the assertion (i) and (1.27) in (ii).
For (1.28) it suffices to prove the following fact.

\proclaim{\noindent Lemma 3.1.}  Let the unit vectors $\omega$ and $\mbox{\boldmath $n$}$ satisfy $\omega\cdot\mbox{\boldmath $n$}\not=0$.  We have
$$\displaystyle
\int_{S(\omega)}\frac{ds(\vartheta)}{(\mbox{\boldmath $n$}\cdot(\omega+i\vartheta))^3}
=\pi(3(\mbox{\boldmath $n$}\cdot\omega)^2-1).
\tag {3.1}
$$

\endproclaim
{\it\noindent Proof.}  The right-hand side on (3.1) is invariant with respect to the change $\omega\rightarrow-\omega$, it is easy to see that 
the case $\omega\cdot\mbox{\boldmath $n$}<0$ can be derived from the result in the case $\omega\cdot\mbox{\boldmath $n$}>0$.
Thus, hereafter we show the validity of (3.1) only for this case.

If $\mbox{\boldmath $n$}\cdot\omega=1$, then $\omega=\mbox{\boldmath $n$}$.  Thus $S(\omega)=S(\mbox{\boldmath $n$})$.  Then
for all $\vartheta\in S(\omega)$ we have $\mbox{\boldmath $n$}\cdot(\omega+i\vartheta)=1$.  This yields
$$\displaystyle
\int_{S(\omega)}\frac{ds(\vartheta)}{(\mbox{\boldmath $n$}\cdot(\omega+i\vartheta))^3}
=2\pi.
$$
Thus the problem is the case when $\mbox{\boldmath $n$}\cdot\omega\not=1$.  Choose an orthogonal $3\times 3$-matrix $A$ such that
$A^T\omega=\mbox{\boldmath $e$}_3$.  Introduce the change of variables $\vartheta=A\vartheta'$.
We have $\vartheta\in S(\omega)$ if and only if $\vartheta'\in S(\mbox{\boldmath $e$}_3)$
and
$$\begin{array}{ll}
\displaystyle
\mbox{\boldmath $n$}\cdot(\omega+iA\vartheta')
&
\displaystyle
=\mbox{\boldmath $n$}'\cdot(\mbox{\boldmath $e$}_3+i\vartheta'),
\end{array}
$$
where $\mbox{\boldmath $n$}'=A^T\mbox{\boldmath $n$}\in S^{2}$.

Here we introduce the polar coordinates for $\vartheta'\in S(\mbox{\boldmath $e$}_3)$:
$$\begin{array}{ll}
\displaystyle
\vartheta'=(\cos\,\vartheta,\sin\varphi, 0)^T, & \varphi\in\,[0,\,2\pi[.
\end{array}
$$
Then, we have
$$\begin{array}{ll}
\displaystyle
I\equiv\int_{S(\omega)}\frac{ds(\vartheta)}{(\mbox{\boldmath $n$}\cdot(\omega+i\vartheta))^3}
&
\displaystyle
=\int_0^{2\pi}\frac{d\varphi}
{(\mbox{\boldmath $n$}'\cdot(i\cos\,\varphi,i\sin\,\varphi,1)^T)^3}
\\
\\
\displaystyle
&
\displaystyle
=-\frac{1}{i}\int_0^{2\pi}\frac{d\varphi}
{(\mbox{\boldmath $n$}'\cdot(\cos\,\varphi,\sin\,\varphi,-i)^T)^3}
\\
\\
\displaystyle
&
\displaystyle
=i\int_0^{2\pi}\frac{d\varphi}
{(a\cos\,\varphi+b\sin\,\varphi-ic)^3},
\end{array}
\tag {3.2}
$$
where $\mbox{\boldmath $n$}'=(a,b,c)^T$.  The numbers $a, b, c$ satisfy $a^2+b^2+c^2=1$ and $0<c<1$
since we have $c=\mbox{\boldmath $n$}'\cdot\mbox{\boldmath $e$}_3=\mbox{\boldmath $n$}\cdot\omega$.
Thus $a^2+b^2\not=0$.  To compute the integral on the right-hand side of (3.2) we make use of the residue calculus.

The change of variables $z=e^{i\varphi}$ gives
$$\begin{array}{l}
\displaystyle
\,\,\,\,\,\,
a\cos\,\varphi+b\sin\,\varphi-ic
\\
\\
\displaystyle
=\frac{1}{2}
\left\{a\left(z+\frac{1}{z}\right)+\frac{b}{i}\left(z-\frac{1}{z}\right)-2ic\right\}
\\
\\
\displaystyle
=\frac{1}{2z}
\left\{(a-ib)z^2-2icz+(a+ib)\right\}
\\
\\
\displaystyle
=\frac{a-ib}{2z}
\left\{\left(z-\frac{ic}{a-ib}\right)^2-\left(\frac{i}{a-ib}\right)^2,\right\}
\\
\\
\displaystyle
=\frac{a-ib}{2z}(z-\alpha)(z-\beta),
\end{array}
\tag {3.3}
$$
where
$$\begin{array}{ll}\displaystyle
\alpha=\frac{i(c+1)}{a-ib}, & 
\displaystyle
\beta=\frac{i(c-1)}{a-ib}.
\end{array}
$$
Since $1-c<1+c$ and $a\cos\,\varphi+b\sin\,\varphi-ic\not=0$ for $z=e^{i\varphi}$, we have
$\vert\beta\vert<1<\vert\alpha\vert$.

Substituting (3.3) into (3.2) and using $d\varphi=\frac{dz}{iz}$, we have
$$\begin{array}{ll}
\displaystyle
I
&
\displaystyle
=i\int_{\vert z\vert=1}\frac{2^3}{(a-ib)^3}\cdot\frac{z^3}{(z-\alpha)^3(z-\beta)^3}\cdot\frac{dz}{iz}
\\
\\
\displaystyle
&
\displaystyle
=\left(\frac{2}{a-ib}\right)^3\int_{\vert z\vert=1}\,\frac{z^2 dz}{(z-\alpha)^3(z-\beta)^3}.
\end{array}
\tag {3.4}
$$
The residue calculus yields
$$\begin{array}{ll}
\displaystyle
\int_{\vert z\vert=1}\,\frac{z^2 dz}{(z-\alpha)^3(z-\beta)^3}
&
\displaystyle
=2\pi i\,\text{Res}_{z=\beta}\,\left(\frac{z^2}{(z-\alpha)^3(z-\beta)^3}\right)
\\
\\
\displaystyle
&
\displaystyle
=2\pi i\cdot\frac{1}{2}\frac{d^2}{dz^2}\left(\frac{z^2}{(z-\alpha)^3}\right)\vert_{z=\beta}
\\
\\
\displaystyle
&
\displaystyle
=2\pi i\cdot\frac{\alpha^2+4\alpha\beta+\beta^2}{(\beta-\alpha)^5}.
\end{array}
\tag {3.5}
$$
Here we have the expression
$$\displaystyle
\alpha-\beta=\frac{2i}{a-ib}
$$
and
$$\begin{array}{ll}
\displaystyle
\alpha^2+4\alpha\beta+\beta^2
&
\displaystyle
=-\frac{(c+1)^2+4(c^2-1)+(c-1)^2}{(a-ib)^2}
\\
\\
\displaystyle
&
\displaystyle
=-\frac{2(3c^2-1)}{(a-ib)^2}.
\end{array}
$$
Thus from (3.4) and (3.5) we obtain 
$$\begin{array}{ll}
\displaystyle
I
&
\displaystyle
=-2\pi\left(\frac{a-ib}{2}\right)^2(\alpha^2+4\alpha\beta+\beta^2)
\\
\\
\displaystyle
&
\displaystyle
=\pi(3c^2-1).
\end{array}
$$
This completes the proof of (3.1).

\noindent
$\Box$

\subsection{On Corollary 1.5}

Let us explain the uniqueness of the solution of the quintic equation (1.30) in $]\frac{1}{\sqrt{3}},\,1]$.

From (1.27), (1.28) and (1.29) we have
$$\displaystyle
\frac{\displaystyle\left\vert\int_{S(\omega)}\,I(\omega,\vartheta)
\,ds(\vartheta)\right\vert}{\pi\,\max_{\vartheta\in S(\omega)}\vert I(\omega,\vartheta)\vert}
=(\mbox{\boldmath $n$}\cdot\omega)^3(3(\mbox{\boldmath $n$}\cdot\omega)^2-1)
$$
and thus
$$\displaystyle
0<
\frac{\displaystyle\left\vert\int_{S(\omega)}\,I(\omega,\vartheta)
\,ds(\vartheta)\right\vert}{\pi\,\max_{\vartheta\in S(\omega)}\vert I(\omega,\vartheta)\vert}
\le 2.
$$
Since $]\,\frac{1}{\sqrt{3}},\,\,1]\ni\gamma\longmapsto\gamma^3(3\gamma^2-1)\in\,]0,\,2]$ is bijective,
the solution of quintic equation (1.30) in $]\frac{1}{\sqrt{3}},\,1]$ is unique and its solution is just
$\gamma=\mbox{\boldmath $n$}\cdot\omega$.

The formulae (1.31) and (1.32) are derived as follows.  A combination of (1.26) and (1.28) yields
$$\displaystyle
(\mbox{\boldmath $n$}\cdot\omega+i\mbox{\boldmath $n$}\cdot\vartheta)^3
=T(\omega,\vartheta).
$$
By expanding the left-hand side, we obtain immediately the desired formulae.

\section{Application to an inverse obstacle problem}

As pointed out in \cite{Ik} the enclosure method developed here can be applied also to an inverse obstacle problem in three dimensions governed by 
the equation
$$\begin{array}{ll}
\displaystyle
\Delta u+k^2 n(x)u=0, & x\in\Omega,
\end{array}
\tag {4.1}
$$
where $k$ is a fixed positive number.  We assume that $\partial\Omega\in C^{\infty}$, for simplicity.
Both $u$ and $n$ can be complex-valued functions.

In this section we assume that $n(x)$ takes the form $n(x)=1+F(x)$, $x\in\Omega$, where
$F=F_{\rho,D}(x)$ is given by (1.7).
We assume that $\rho\in L^{\infty}(D)$ instead of $\rho\in L^2(D)$ and that $u\in H^2(\Omega)$ is an arbitrary non trivial solution of (4.1)
at this stage.
We never specify the boundary condition of $u$ on $\partial\Omega$. 
By the Sobolev imbedding theorem \cite{G} one may assume that $u\in C^{0,\alpha}(\overline\Omega)$ with $0<\alpha<1$.

In this section we consider

{\bf\noindent Problem 2.} Extract information about the singularity of $D$ from the Cauchy data of $u$ on $\partial\Omega$.

We encounter this type of problem, for example, $u$ is given by the restriction to $\Omega$ of the total wave
defined in the whole space and generated by a point source located outside of $\Omega$ or a single plane wave coming from infinity.
The surface where the measurements are taken is given by $\partial\Omega$ which encloses the penetrable obstacle $D$
with a different reflection index $1+\rho$, $\rho\not\equiv 0$.  See \cite{CK} for detailed information about the direct problem itself.
Any way we start with having the Cauchy data of an arbitrary (nontrivial) $H^2(\Omega)$ solution of (4.1).

Using the Cauchy data of $u$ on $\partial\Omega$, we introduce the indicator function
$$\displaystyle
I_{\omega,\vartheta}(\tau)=\int_{\partial\Omega}
\left(\frac{\partial u}{\partial\nu}v-\frac{\partial v}{\partial\nu} u\right)\,dS,
\tag {4.2}
$$
where the function $v=v(x), x\in\Bbb R^3$ is given by
$$\displaystyle
v=e^{x\cdot(\tau\omega+i\sqrt{\tau^2+k^2}\vartheta)},\,\,\tau>0
$$
and $\vartheta\in S(\omega)$.
And also its derivative with respect to $\tau$ is given by the formula
$$\displaystyle
I_{\omega,\vartheta}'(\tau)
=\int_{\partial\Omega}\left(\frac{\partial u}{\partial\nu}\,v_{\tau}-\frac{\partial\,v_{\tau}}{\partial\nu} u\right)\,dS,
\tag {4.3}
$$
where
$$\displaystyle
v_{\tau}=\partial_{\tau}v=\left\{x\cdot\left(\omega+i\frac{\tau}{\sqrt{\tau^2+k^2}}\,\vartheta\,\right)\,\right\}\,v.
$$

As done the proof of Theorem 1.1 integration by parts yields
$$\displaystyle
I_{\omega,\vartheta}(\tau)=-k^2\int_D\rho(x)u(x)v\,dx
$$
and
$$\displaystyle
I_{\omega,\vartheta}'(\tau)=-k^2\int_D\rho(x)u(x)v_{\tau}\,dx.
$$
Thus this can be viewed as the case $\rho(x)$ in Problem 1 is given by $-k^2\rho(x)u(x)$ and $\tilde{\rho}(x)$ in Definition 1.2
by $-k^2\tilde{\rho}(x)u(x)$.

Thus we obtain
\proclaim{\noindent Theorem 4.1.}
Let $\omega$ be regular with respect to $D$ and assume that
$D$ has a conical singularity from direction $\omega$.
Then, we have
$$\displaystyle
\tau^3e^{-\tau h_D(\omega)}e^{-i\sqrt{\tau^2+k^2}p(\omega)\cdot\vartheta}I_{\omega,\vartheta}(\tau)=
-k^2\tilde{\rho}(p(\omega))\,u(p(\omega))\,C_{(p(\omega),\omega)}(\delta,Q,\vartheta)
+O(\tau^{-\alpha})
$$
and
$$\displaystyle
\tau^3e^{-\tau h_D(\omega)}e^{-i\sqrt{\tau^2+k^2}p(\omega)\cdot\vartheta}I_{\omega,\vartheta}'(\tau)=
-k^2\tilde{\rho}(p(\omega))\,u(p(\omega))(h_D(\omega)+ip(\omega)\cdot\vartheta)\,C_{(p(\omega),\omega)}(\delta,Q,\vartheta)
+O(\tau^{-\alpha}).
$$
The remainder $O(\tau^{-\alpha})$ is uniform with respect to $\vartheta\in S(\omega)$.
\endproclaim

Thus under the same assumptions as Theorem 4.1, for each $\vartheta\in S(\omega)$ one can calculate 
$$\displaystyle
I(\omega\,\vartheta)\equiv -k^2\tilde{\rho}(p(\omega))\,u(p(\omega))\,C_{(p(\omega),\omega)}(\delta,Q)
$$
via the formula
$$\displaystyle
I(\omega,\vartheta)
=\lim_{\tau\rightarrow\infty}\tau^3e^{-\tau h_D(\omega)}e^{-i\sqrt{\tau^2+k^2}p(\omega)\cdot\vartheta}
I_{\omega,\vartheta}(\tau)
\tag {4.4}
$$
by using the Cauchy data of $u$ on $\partial\Omega$ if $p(\omega)$ is known.

And also we have
\proclaim{\noindent Theorem 4.2.}
Let $\omega$ be regular with respect to $D$.
Assume that $D$ has a conical singularity from direction $\omega$;
$n(x)-1=F_{\rho,D}(x)$ is active at $p(\omega)$ in the sense of
Definition 1.2 and the value of $u$ at $p(\omega)$ satisfies
$$\displaystyle
u(p(\omega))\not=0.
\tag {4.5}
$$
If the direction $\vartheta\in S(\omega)$ satisfies the condition (1.15), then all the formulae (1.16), (1.17)
and (1.18) for the indicator function defined by (4.2) together with its derivative (4.3) are valid.

\endproclaim
Note that the assumption (4.5) ensures $u\not\equiv 0$. See Appendix for an example of $u$ satisfying (4.5).

The following corollaries corresponds to Corollaries 1.1 and 1.2.

\proclaim{\noindent Corollary 4.1.}
Let $\omega$ be regular with respect to $D$.
Under the same assumptions as those in Theorem 4.2 the point $p(\omega)$ is uniquely determined by
the Cauchy data of $u$ on $\partial\Omega$.

\endproclaim

\proclaim{\noindent Corollary 4.2.}  Let $u\in H^2(\Omega)$ be a solution of (4.1).
Assume that $D$ is given by the inside of a convex polyhedron and that in a neighbourhood of each vertex $p$ of $D$,  
the $D$ coincides with the inside of a tetrahedron with apex $p$ and that $n-1=F_{\rho, D}$ given by (1.7)
is active at $p$ and the value of $u$ at $p$ satisfies (4.5).
Then, all the formulae (1.16), (1.17) and (1.18) for the indicator function defined by (4.2) together with its derivative (4.3)
are valid for all $\omega$ regular with respect to $D$
and $\vartheta\in S(\omega)$.
Besides, the Cauchy data of $u$ on $\partial\Omega$ uniquely determines $D$.
\endproclaim

The following result is an extension of Theorem 4.1 in \cite{Ik} to three dimensional case.

\proclaim{\noindent Corollary 4.3.}  Let $u\in H^2(\Omega)$ be a solution of (4.1).
Let $\omega\in S^2$ be regular with respect to $D$.  Assume that: $D$ has a circular cone singularity
at $p=p(\omega)$; $n(x)-1=F_{\rho,D}(x)$ is active at $p(\omega)$ in the sense of
Definition 1.2 and the value of $u$ at $p(\omega)$ satisfies (4.5).
Choose two linearly independent vectors $\vartheta=\vartheta_1$ and $\vartheta_2$ in $S(\omega)$.
Then, the point $p(\omega)$ itself and thus $h_D(\omega)=p(\omega)\cdot\omega$ can be extracted from the Cauchy data of $u$ on $\partial\Omega$
by using the formula
$$\displaystyle
p(\omega)\cdot\omega+i\,p(\omega)\cdot\vartheta_j
=\lim_{\tau\rightarrow\infty}
\frac{I_{\omega,\vartheta_j}'(\tau)}{I_{\omega,\vartheta_j}(\tau)},\,\,\,j=1,2.
\tag {4.6}
$$
\endproclaim

By virtue of the formula (1.24), the function $I(\omega,\,\cdot\,)$ has the expression
$$\displaystyle
I(\omega,\vartheta)=-6k^2\,\tilde{\rho}(p(\omega))u(p(\omega))\,V(\theta)(\mbox{\boldmath $n$}\cdot(\omega+i\vartheta))^{-3}.
\tag {4.7}
$$
Simillarly to Corollary 1.4 formula (4.7) yields immediately the following results.

\proclaim{\noindent Corollary 4.4.}  Let $u\in H^2(\Omega)$ be a solution of (4.1).
Let $\omega\in S^2$ be regular with respect to $D$.  Assume that: $D$ has a circular cone singularity
at $p(\omega)$ such as $D\cap B_{\epsilon}(p(\omega))=V_{p(\omega)}(-\mbox{\boldmath $n$},\theta)\cap B_{\epsilon}(p(\omega))$ with a $\epsilon>0$.

\noindent
(i)  Assume that $n(x)-1=F_{\rho,D}(x)$ is active at $p(\omega)$ in the sense of
Definition 1.2 and the value of $u$ at $p(\omega)$ satisfies (4.5).
The vector $\omega$ coincides with $\mbox{\boldmath $n$}$ if and only if the function
$I(\omega,\,\cdot\,)$ is a constant function.

\noindent
(ii)  The vector $\mbox{\boldmath $n$}$ and $\theta$ of $V_{p(\omega)}(-\mbox{\boldmath $n$},\theta)$
and $\tilde{\rho}(p(\omega))\,u(p(\omega))$ satisfies the following two equations:
$$\displaystyle
6k^2\,\vert\tilde{\rho}(p(\omega))\,u(p(\omega))\vert\,V(\theta)=(\mbox{\boldmath $n$}\cdot\omega)^3
\max_{\vartheta\in S(\omega)}\vert I(\omega,\vartheta)\vert;
\tag {4.8}
$$
$$\displaystyle
-6k^2\,\tilde{\rho}(p(\omega))u(p(\omega))
\,V(\theta)\,(3(\mbox{\boldmath $n$}\cdot\omega)^2-1)
=\frac{1}{\pi}\,\int_{S(\omega)}\,I(\omega,\vartheta)
\,ds(\vartheta).
\tag {4.9}
$$

\endproclaim

Using the equations (4.7), (4.8) and (4.9) one gets the following corollary.

\proclaim{\noindent Corollary 4.5.}  Let $u\in H^2(\Omega)$ be a solution of (4.1).
Let $\omega\in S^2$ be regular with respect to $D$.  Assume that: $D$ has a circular cone singularity
at $p(\omega)$ such as $D\cap B_{\epsilon}(p(\omega))=V_{p(\omega)}(-\mbox{\boldmath $n$},\theta)\cap B_{\epsilon}(p(\omega))$ with a $\epsilon>0$.
Assume that $n(x)-1=F_{\rho,D}(x)$ is active at $p(\omega)$ in the sense of
Definition 1.2 and the value of $u$ at $p(\omega)$ satisfies (4.5).
Assume also that $\omega\approx\mbox{\boldmath $n$}$ in the sense that (1.29) holds.
Then, we have the completely same statement and formulae as those of Corolloary 1.5.
\endproclaim

Note that under the same assumptions as Corollary 4.5, one can finally calculate the quantity
$$\displaystyle
\tilde{\rho}(p(\omega))\,u(p(\omega))\,V(\theta)
\tag {4.10}
$$
and $\mbox{\boldmath $n$}$ from the Cauchy data of $u$ on $\partial\Omega$.
Since the steps for the calculation are similar to the steps presented in Subsection 1.2 for the inverse souce problem, we omit its description.

However, it should be noted that, in addition, if $\tilde{\rho}(p(\omega))$ is known to be a {\it real number},
then one can recover the phase of the complex number $u(p(\omega))$ modulo $2\pi n$, $n=0,\pm 1,\pm 2,\cdots$ from the computed value (4.10).

{\bf\noindent Remark 4.1.}
One can apply the result in \cite{IkC} to the computation of the value $u(p(\omega))$ itself.
For simplicity we assume that $\Omega$ is convex like a case when $\Omega=B_R(x_0)$ centered at a point $x_0$ with a large radius $R$.
From formula (4.6) we know the position of $p(\omega)$ and thus the domain $\Omega\cap\{x\in\Bbb R^3\,\vert\,x\cdot\omega>x\cdot p(\omega)\}$.
Because of the continuity of $u$ on $\overline\Omega$, one has, for a sufficiently small $\epsilon>0$
$$\displaystyle
u(p(\omega))\approx u(p(\omega)+\epsilon\,\omega).
$$
Since the point $p(\omega)+\epsilon\,\omega\in\Omega\cap\{x\in\Bbb R^3\,\vert\,x\cdot\omega>x\cdot p(\omega)\}$
and therein $u\in H^2$ satisfies the Helmholtz equation $\Delta u+k^2u=0$, one can calculate the value
$u(p(\omega)+\epsilon\,\omega)$ itself from the Cauchy data of $u$ on $\partial\Omega\cap\{x\in\Bbb R^3\,\vert\,x\cdot\omega>x\cdot p(\omega)\}$
by using Theorem 1 in \cite{IkC}.

\section{Final remark}

All the results in this paper can be extended also to the case when the governing equation of the background 
medium is given by a Helmholtz equation with a known coefficient $n_0(x)$.  It means that
if one considers, instead of (1.6) and (4.1) the equations
$$\begin{array}{ll}
\displaystyle
\Delta u+k^2n_0(x)u=F_{\rho,D}(x), & x\in\Omega
\end{array}
$$
and
$$\begin{array}{ll}
\displaystyle
\Delta u+k^2(n_0(x)+F_{\rho,D}(x))u=0, & x\in\Omega,
\end{array}
$$
respectively, then one could obtain all the corresponding results.

\section{Appendix.  On condition (4.5)}

As suggested in \cite{Ik} the condition (4.5) can be satisfied if $k$ is sufficiently small
under the situation when $u$ is given by the restriction onto $\Omega$ of the total 
field $U$ in the whole space scattering problem generated by, for example, the point source located at a point $z$ in $\Bbb R^3\setminus\overline\Omega$.
The $U$ has the expression $U=\Phi(x,z)+w_z(x)$, where
$$\begin{array}{ll}
\displaystyle
\Phi(x,z)=\frac{1}{4\pi}\frac{e^{ik\vert x-z\vert}}{\vert x-z\vert}, & x\in\Bbb R^3\setminus\{z\}
\end{array}
$$
and $w_z\in H^2_{\text{local}}(\Bbb R^3)$ is the unique solution of the inhomogeneous Helmholtz equation
$$\begin{array}{ll}
\displaystyle
\Delta w_z+k^2w_z+k^2F(x)(w_z+\Phi(x,z))=0, & x\in\Bbb R^3
\end{array}
$$
with the outgoing Sommerfeld radiation condition
$$\displaystyle
\lim_{r\rightarrow\infty}r\left(\frac{\partial}{\partial r}w_z(x)-ik w_z(x)\right)=0,
$$
where $r=\vert x\vert$ and $F=F_{\rho,D}$ is given by (1.7).  See \cite{CK} for the solvabilty.

Here we claim
\proclaim{\noindent Proposition A.}
Let $0<R_1<R_2$ satisfy $D\subset B_{R_2}(z)\setminus\overline B_{R_1}(z)$.
Let $M>0$ and $R>0$ satisfy $\vert D\vert\le M$ and $\Vert\rho\Vert_{L^{\infty}(D)}\le R$, respectively.

If $k$ satisfies the system of inequalities
$$\displaystyle
C\equiv \frac{3k^2R}{2}\left(\frac{M}{4\pi}\right)^{2/3}<1
\tag {A.1}
$$
and
$$\displaystyle
\frac{C}{1-C}<\frac{R_1}{R_2},
\tag {A.2}
$$
then, for all $x\in\overline D$ we have
$$\displaystyle
\vert U(x)\vert\ge
\frac{1}{4\pi}\left(\frac{1}{R_2}-\frac{C}{1-C}\frac{1}{R_1}\,\right).
\tag {A.3}
$$

\endproclaim

{\it\noindent Proof.}  Note that $w\in C^{0,\alpha}(\overline\Omega)$ with $0<\alpha<1$ by the Sobolev imbedding theorem.
It is well known that the function $w_z$ satisfies the Lippman-Schwinger equation
$$\begin{array}{ll}
\displaystyle
w_z(x)
&
\displaystyle
=k^2\int_{D}\Phi(x,y)\rho(y)w_z(y)\,dy+k^2\int_{D}\Phi(x,y)\Phi(y,z)\rho(y)\,dy
\end{array}
$$
and thus, for all $x\in\overline{D}$ we have
$$\displaystyle
\vert w_z(x)\vert
\le
\frac{k^2 R}{4\pi}
\left(\Vert w_z\Vert_{L^{\infty}(D)}
+\frac{1}{4\pi\,R_1}\right)\,
\int_D\frac{dy}{\vert x-y\vert}.
\tag {A.4}
$$
Let $\epsilon>0$.  We have
$$\begin{array}{ll}
\displaystyle
\int_D\frac{dy}{\vert x-y\vert}
&
\displaystyle
=\int_{D\cap B_{\epsilon}(x)}\,\frac{dy}{\vert x-y\vert}+\int_{D\setminus B_{\epsilon}(x)}\,\frac{dy}{\vert x-y\vert}
\\
\\
\displaystyle
&
\displaystyle
\le
\int_{B_{\epsilon}(x)}\,\frac{dy}{\vert x-y\vert}+\frac{\vert D\vert}{\epsilon}
\\
\\
\displaystyle
&
\displaystyle
\le
2\pi\epsilon^2+\frac{M}{\epsilon}.
\end{array}
$$
Choose $\epsilon$ in such a way that this right-hand side becomes minimum, that is,
$$\displaystyle
\epsilon=\left(\frac{M}{4\pi}\right)^{1/3}.
$$
Then one gets
$$\begin{array}{ll}
\displaystyle
\int_D\frac{dy}{\vert x-y\vert}
&
\displaystyle
\le
6\pi
\left(\frac{M}{4\pi}\right)^{2/3}.
\end{array}
$$
Thus this together with (A.4) yields
$$\displaystyle
\left(1-C\,\right)\Vert w_z\Vert_{L^{\infty}(D)}
\le
\frac{C}{4\pi\,R_1}.
$$
This together with the estimate
$$\displaystyle
\vert U(x)\vert\ge \frac{1}{4\pi\,R_2}-\Vert w_z\Vert_{L^{\infty}(D)}
$$
yields the desired estimate (A.3) under the assumptions (A.1) and (A.2).

\noindent
$\Box$

Note that since $R_2>R_1$, the set of inequalities (A.1) and (A.2) are equivalent to the single inequality
$$\displaystyle
C<\frac{R_1}{R_1+R_2}.
\tag {A.5}
$$
Thus we choose $k^2$ sufficiently small in the sense of (A.5) we have, for all $x\in\overline D$
$$\displaystyle
\vert u(x)\vert \ge\frac{1}{4\pi}\left(\frac{1}{R_2}-\frac{C}{1-C}\frac{1}{R_1}\,\right)>0.
$$
Thus the condition (4.5) for $u=U\vert_{\Omega}$ is satisfied.  The choice of $k$ depends only on the a-priori information
about $D$ and $\rho$
described by $R_1$, $R_2$, $M$ and $R$.

$$\quad$$

\centerline{{\bf Acknowledgments}}

This research was partially supported by Grant-in-Aid for Scientific Research
(C) (No. 17K05331) and (B) (N. 18H01126) of Japan Society for the Promotion of Science.

$$\quad$$

\end{document}